\documentclass{amsart}

\newtheorem{teo}{Theorem}
\newtheorem{lem}{Lemma}
\newtheorem{pro}{Proposition}
\newtheorem{cor}{Corollary}

\newtheorem{defi}{Definition}

\theoremstyle{definition}

\newtheorem{examp}{Example}

\def\sp{\hspace{0.2cm}}
\def\qed{\begin{flushright} $\square$ \end{flushright}}
\def\qee{\begin{flushright} $\Diamond$ \end{flushright}}
\def\ov{\overline}

\begin{document}

\title[A Thermodynamic Formalism for density matrices in Q. I.]
{ A Thermodynamic Formalism for density matrices in Quantum
Information}


\author{A. Baraviera}
\address{}
\curraddr{}
\email{}
\thanks{}

\author{C. F. Lardizabal}
\address{}
\curraddr{}
\email{}
\thanks{Supported in part by CAPES and CNPq}

\author{A. O. Lopes}
\address{}
\curraddr{}
\email{}
\thanks{}

\author{M. Terra Cunha}
\address{}
\curraddr{}
\email{}
\thanks{}

\begin{abstract}

We consider new concepts of entropy and pressure for stationary
systems acting on density matrices  which generalize the usual
ones in Ergodic Theory. Part of our work is to justify why the
definitions and results we describe here are natural
generalizations of the classical concepts of Thermodynamic
Formalism (in the sense of R. Bowen, Y. Sinai and D. Ruelle). It
is well-known that the concept of density operator should replace
the concept of measure for the cases in which we consider a
quantum formalism.

We consider the
operator $\Lambda$ acting on the space of density matrices $\mathcal{M}_N$ over a finite $N$-dimensional complex Hilbert space
$$
\Lambda(\rho):=\sum_{i=1}^k tr(W_i\rho W_i^*)\frac{V_i\rho V_i^*}{tr(V_i\rho V_i^*)},
$$
where $W_i$ and $V_i$, $i=1,2,\dots, k$ are linear operators in this Hilbert
space. In some sense this operator is a version of an Iterated Function System (IFS). Namely, the $V_i\,(.)\,V_i^*=:F_i(.)$, $i=1,2,\dots,k$, play the role of the inverse branches (i.e., the dynamics on the configuration space of density matrices) and the $W_i$ play the role of the weights one can consider on the IFS. In this way a family $W:=\{W_i\}_{i=1,\dots, k}$ determines a Quantum Iterated Function System (QIFS).

We  also present some estimates related to the Holevo bound.
\vspace{0.2cm}

\end{abstract}

\maketitle

\normalsize{

\bigskip

{Paper to appear in Applied Mathematics Research Express (2010)}

\bigskip

\section{Introduction}}

In this work we investigate a generalization of the
classical Thermodynamic Formalism (in the sense of Bowen, Sinai and
Ruelle) for the setting of density matrices. We consider the
operator $\Lambda$ acting on the space of density matrices $\mathcal{M}_N$ over a finite $N$-dimensional complex Hilbert space
\begin{equation}
\Lambda(\rho):=\sum_{i=1}^k tr(W_i\rho W_i^*)\frac{V_i\rho V_i^*}{tr(V_i\rho V_i^*)},
\end{equation}
where $W_i$ and $V_i$, $i=1,2,\dots, k$ are linear operators in this Hilbert
space. Note that $\Lambda$ is not a linear operator. This operator can be seen as a version of an Iterated Function System (IFS). Namely, the $V_i\,(.)\,V_i^*=:F_i(.)$, $i=1,2,\dots,k$, play the role of the inverse branches (i.e., the dynamics on the configuration space of density matrices $\rho$) and the $W_i$ play the role of the weights one can consider on the IFS.
We suppose that for all $\rho$ we have that $\sum_{i=1}^k tr(W_i\rho W_i^*)=1$.  Note that such trace preserving condition, for any normalized operator $\rho$ (that is, with $tr(\rho)= 1$), is equivalent to the explicit condition $\sum_i W_i^* W_i = I$. We say that $\Lambda$ is a normalized operator.

A  family $W:=\{W_i\}_{i=1,\dots, k}$ determines a Quantum Iterated Function System (QIFS) $\mathcal{F}_{W}$,
$$\mathcal{F}_W=\{\mathcal{M}_N,F_i,W_i\}_{i=1,\dots, k}$$
Basic references on QIFS are \cite{jadcyk} and \cite{lozinski}. We want to consider a new concept of entropy for stationary systems acting on density matrices which generalizes the usual one in
Ergodic Theory. In our setting the $V_i$, $i=1,2,\dots,k$ are fixed
(i.e. the dynamics of the inverse branches is fixed in the beginning) and we
consider the different families $W_i$, $i=1,2,\dots,k$, (also with
the attached corresponding eigendensity matrix $\rho_W$) as possible
Jacobians of stationary probabilities.

Given a normalized family $W_i$, $i=1,2,\dots ,k$, a natural definition of entropy is given by
\begin{equation}
h_V(W)=-\sum_{i=1}^k \frac{tr(W_i\rho_W W_i^*)}{tr(V_i\rho_W V_i^*)}\sum_{j=1}^k tr\Big(W_j V_i\rho_W V_i^*
W_j^*\Big)\log{\Big(\frac{tr(W_j V_i\rho_W V_i^* W_j^*)}{tr(V_i\rho_W V_i^*)}\Big)}
\end{equation}
where $\rho_W$ denotes the barycenter of the unique invariant,
attractive measure for the Markov operator $\mathcal{V}$ associated
to $\mathcal{F}_W$. We show that this generalizes the entropy of a
Markov System.

We also want to present a concept of pressure for stationary systems acting on density matrices  which generalizes the usual one in Ergodic Theory. In addition to the dynamics obtained by the $V_i$, which are fixed, a family of potentials $H_i$,  $i=1,2,\dots k$ induces a kind of Ruelle operator given by
\begin{equation}\label{operador_usa1a}
\mathcal{L}_H(\rho):=\sum_{i=1}^k tr(H_i\rho H_i^*)V_i\rho V_i^*
\end{equation}

We show that such operator admits an eigenvalue $\beta$ and an
associated eigenstate $\rho_\beta$, that is, one satisfying
$\mathcal{L}_H(\rho_\beta)=\beta\, \rho_{\beta}$.

The natural generalization of the concept of pressure for a family
$H_i$, $i=1,2,\dots ,k$ is the problem of maximizing, on
the possible normalized families $W_i$, $i=1,2,\dots k$, the
expression
\begin{equation}
h_V(W)+ \sum_{j=1}^k \log \Big( tr(H_j\rho_\beta
H_j^*)tr(V_j\rho_\beta V_j^*)\Big) tr(W_j \rho_W  W_j^*)
\end{equation}
We show a relation between the eigendensity matrix $\rho_\beta$ for the Ruelle
operator and the set of $W_i$, $i=1,2,\dots k$, which maximizes
pressure. In the particular case that each of the $V_i$ is unitary, $i=1,2,\dots k$, the maximum value is $\log \beta$.

Our work is inspired by the results presented in \cite{lozinski} and
\cite{wsbook}. We would like to thank these authors for supplying us with the
corresponding references.

It is well-known that  completely positive mappings (operators)
acting on density matrices are of great importance in Quantum
Computing. These operators can be written  in the
Stinespring-Kraus form (see section \ref{sec_cc}). Also a nice
exposition on the interplay of Ergodic Theory and Quantum
Information is presented in \cite{be}.

The initial part of our work aims to present some of the
definitions and concepts that are not very well-known (at least
for the general audience of people in Dynamical Systems), in a
systematic way. We present the main basic definitions which are
necessary to understand the theory. However, we do not have the
intention of exhausting what is already known. We believe that the
theoretical results presented here can be useful as a general tool
to understand problems in Quantum Computing.

Several examples are presented in the text.  We believe  that this
will help the reader to understand some of the main issues of the theory. In order to simplify the notation we will present most of our results for the case of matrices of order 2.

In sections \ref{sec_basic} and  \ref{sec_examples} we present some basic definitions, examples and we show some preliminary relations of our setting to the classical Thermodynamic Formalism. In section \ref{sec_eigenvalues} we present an eigenvalue problem for non-normalized Ruelle operators  which will be required later. Some properties and concepts about density matrices and Ruelle operators are presented in sections \ref{sec_densities} and \ref{sec_somelemmas}. Sections \ref{sec_entr_ifs} and \ref{sec_calc_entr}  are dedicated to the introduction of some different kinds of entropy that were already known but do not have a stationary character. In section \ref{sec_novaentr} we introduce the concept of stationary entropy for {\it measures} defined on the set of density matrices. In section \ref{camarkov} we compare this definition with the usual one for Markov Chains.
Section \ref{sec_cc} is dedicated to motivate the interest on pressure and the capacity-cost function. Section \ref{sec_analysis1}, \ref{sec_analysis2}, \ref{sec_analysis3} and \ref{sec_analysis4} are dedicated to the presentation of our main results on pressure, important inequalities, examples and its relation with the classical theory of Thermodynamic Formalism.

In \cite{BLLT1} we present a general exposition (describing the setting we consider here) where we omit proofs, but provide many examples. We believe that paper will help to complement the present paper for the reader which is a newcomer in the area. We also present there some basic results concerning the discrete Wigner measure.

In \cite{BLLT2} we propose a different concept of entropy which is also a generalization of the classical one. We also describe some properties of the Quantum Stochastic Process associated to the Quantum Iterated Function System.

This work is part of the thesis dissertation of C. F. Lardizabal in Prog. Pos-Grad. Mat. UFRGS (Brazil).

\bigskip

\section{Basic definitions}\label{sec_basic}

Let $M_N(\mathbb{C})$ the set of complex matrices  of order $n$. If $\rho\in M_N(\mathbb{C})$ then $\rho^*$ denotes the transpose conjugate of $\rho$. A state (or vector) in $\mathbb{C}^n$ will be denoted by $\psi$ or $\vert \psi\rangle$, and the associated projection will be written $\vert \psi\rangle\langle\psi\vert$. Define
$$\mathcal{H}_N:=\{\rho\in M_N(\mathbb{C}):\rho^*=\rho\}$$
$$\mathcal{PH}_N:=\{\rho\in \mathcal{H}_N:\langle\rho \psi,\psi\rangle\geq 0 , \forall \psi\in\mathbb{C}^N\}$$
$$\mathcal{M}_N:=\{\rho\in \mathcal{PH}_N: tr(\rho)=1\}$$
$$\mathcal{P}_N:=\{\rho\in \mathcal{H}_N: \rho=\vert\psi\rangle\langle\psi\vert,\psi\in\mathbb{C}^N, \langle\psi\vert\psi\rangle=1\},$$
the space of hermitian, positive, density operators and pure states, respectively. Density operators are also called mixed states. If a quantum system can be in one of the states $\{\psi_1,\dots, \psi_k\}$ then a mixed state $\rho$ will be written as
\begin{equation}
\rho=\sum_{i=1}^k p_i\vert \psi_i\rangle\langle \psi_i\vert
\end{equation}
where the $p_i$ are positive numbers with $\sum_i p_i=1$.

\bigskip

\begin{defi}
Let $F_i:\mathcal{M}_N\to\mathcal{M}_N$, $p_i:\mathcal{M}_N\to[0,1]$, $i=1,\dots ,k$ and such that $\sum_i p_i(\rho)=1$. We call
\begin{equation}\label{qifs}
\mathcal{F}_N=\{\mathcal{M}_N, F_i, p_i:i=1,\dots, k\}
\end{equation}
a {\bf Quantum Iterated Function System} (QIFS).
\end{defi}

\begin{defi}
A QIFS is {\bf homogeneous} if $p_i$ and $F_ip_i$ are affine mappings, $i=1,\dots , k$.
\end{defi}

Suppose that the QIFS considered is such that there are $V_i$ and $W_i$ linear maps, $i=1,\dots, k$, with $\sum_{i=1}^k W_i^*W_i=I$ such that
\begin{equation}
F_i(\rho)=\frac{V_i\rho V_i^*}{tr(V_i\rho V_i^*)}
\end{equation}
and
\begin{equation}
p_i(\rho)=tr(W_i\rho W_i^*)
\end{equation}
Then we have that a QIFS is homogeneous if $V_i$=$W_i$,
$i=1,\dots,k$. Now we can define a Markov operator
$\mathcal{V}:\mathcal{M}^1(\mathcal{M}_N)\to\mathcal{M}^1(\mathcal{M}_N)$,
\begin{equation}
(\mathcal{V}\mu)(B)=\sum_{i=1}^k\int_{F_i^{-1}(B)}p_i(\rho)d\mu(\rho),
\end{equation}
where $\mathcal{M}^1(\mathcal{M}_N)$ denotes the space of probability measures over $\mathcal{M}_N$. We also define $\Lambda: \mathcal{M}_N\to\mathcal{M}_N$,
\begin{equation}
\Lambda(\rho):=\sum_{i=1}^k p_i(\rho)F_i(\rho)
\end{equation}
The operator defined above has no counterpart in the classical Thermodynamic Formalism. We will also consider the operator defined on the space of density matrices $\rho$,
\begin{equation}
\mathcal{L}(\rho)=\sum_{i=k}^k q_i(\rho)V_i\rho V_i^*.
\end{equation}
If for all $\rho$ we have $\sum_{i=k}^k q_i(\rho)=1$, we say the operator is normalized. We are also interested in the non-normalized case. If the QIFS is homogeneous, then
\begin{equation}\label{lambda_um}
\Lambda(\rho)=\sum_i V_i\rho V_i^*
\end{equation}

\begin{teo}\cite{wsbook} A mixed state $\hat{\rho}$ is $\Lambda$-invariant if and only if
\begin{equation}\label{bari}
\hat{\rho}=\int_{\mathcal{M}_N} \rho d\mu(\rho),
\end{equation}
for some $\mathcal{V}$-invariant measure $\mu$.
\end{teo}
We recall the definition of the integral above in section \ref{v_bari}.

\bigskip

In order to define hyperbolic QIFS, one has to define a distance on the space of mixed states. For instance, we could choose one of the following:
\begin{equation}
D_1(\rho_1,\rho_2)=\sqrt{tr[(\rho_1-\rho_2)^2]}
\end{equation}
\begin{equation}
D_2(\rho_1,\rho_2)=tr\sqrt{(\rho_1-\rho_2)^2}
\end{equation}
\begin{equation}
D_3(\rho_1,\rho_2)=\sqrt{2\{1-tr[(\rho_1^{1/2}\rho_2\rho_1^{1/2})^{1/2}]\}},
\end{equation}
the Hilbert-Schmidt, trace, and Bures distances, respectively. Such metrics generate the same topology on $\mathcal{M}_N$. Considering the space of mixed states with one of those metrics we can use a definition of hyperbolicity similar to the one used for IFS. That is, we say a QIFS is {\bf hyperbolic} if the quantum maps $F_i$ are contractions with respect to one of the distances on $\mathcal{M}_N$ and if the maps $p_i$ are H\"older-continuous and positive, see for instance \cite{lozinski}.

\bigskip

\begin{pro}\label{hh_teo_i} If a QIFS (\ref{qifs}) is homogeneous and hyperbolic then the associated Markov operator admits a unique invariant measure $\mu$. Such invariant measure determines a unique  $\Lambda$-invariant state $\rho\in\mathcal{M}_N$, given by (\ref{bari}).
\end{pro}

See \cite{lozinski}, \cite{wsbook} for the proof.

\section{Examples of QIFS}\label{sec_examples}

\begin{examp}\label{exemp7} $\Omega=\mathcal{M}_N$, $k=2$, $p_1=p_2=1/2$, $G_1(\rho)=U_1\rho U_1^*$, $G_2(\rho)=U_2\rho U_2^*$. The normalized identity matrix $\rho_*=I/N$ is $\Lambda$-invariant, for any choice of unitary $U_1$ and $U_2$. Note that we can write
\begin{equation}
\rho_*=\int_{\mathcal{M}_N}\rho d\mu(\rho)
\end{equation}
where the measure $\mu$, uniformly distributed over
$\mathcal{P}_N$ (the Fubini-Study metric), is
$\mathcal{V}$-invariant.
\end{examp}

\qee

We recall that a mapping $\Lambda$ is {\bf completely positive} (CP) if $\Lambda\otimes I$ is positive for any extension of the Hilbert space considered $\mathcal{H}_N\to\mathcal{H}_N\otimes\mathcal{H}_E$. We know that every CP mapping which is trace-preserving can be represented (in a nonunique way) in the Stinespring-Kraus form
\begin{equation}
\Lambda(\rho)=\sum_{j=1}^k V_j\rho V_j^*, \sp \sum_{j=1}^k V_j^* V_j=1,
\end{equation}
where the $V_i$ are linear operators. Moreover if we have $\sum_{j=1}^k V_j  V_j^*=I$, then $\Lambda(I/N)=I/N$. This is the case if each of the $V_i$ are normal.

\bigskip

We call a unitary trace-preserving CP map a {\bf bistochastic map}. An example of such a mapping is
\begin{equation}
\Lambda_U(\rho)=\sum_{i=1}^k p_i U_i\rho U_i^*,
\end{equation}
where the $U_i$ are unitary operators and $\sum_i p_i=1$. Note that if we write $F_i(\rho)=U_i\rho U_i^*$, then example \ref{exemp7} is part of this class of operators. For such operators we have that $\rho_*$ is an invariant state for  $\Lambda_U$ and also that $\delta_{\rho_*}$ is invariant for the Markov operator $P_U$ induced by this QIFS.

\bigskip
We will present a simple example of the kind of problems we are interested here, namely eigenvalues and eigendensity matrices. Let $\mathcal{H}_N$ be a Hilbert space of dimension $N$.  As before, let $\mathcal{M}_N$ be the space of density operators on $\mathcal{H}_N$. A natural problem is to find fixed points for $\Lambda:\mathcal{M}_N\to\mathcal{M}_N$,
\begin{equation}
\Lambda(\rho)=\sum_{i=1}^k V_i\rho V_i^*
\end{equation}

In order to simplify our notation we fix $N=2$ and  $k=2$. Let
$$V_1=\left(
\begin{array}{cc}
v_1 & v_2\\
v_3 & v_4
\end{array}
\right),
\sp V_2=\left(
\begin{array}{cc}
w_1 & w_2\\
w_3 & w_4
\end{array}
\right),
\sp \rho=\left(
\begin{array}{cc}
\rho_1 & \rho_2\\
\ov{\rho_2} & \rho_4
\end{array}
\right),
$$
 where $V_1$ and $V_2$ are invertible and $\rho$ is a density operator. We would like to find $\rho$ such that
\begin{equation}\label{equm}
V_1\rho V_1^*+V_2\rho V_2^*=\rho.
\end{equation}

\begin{examp}\label{vex1} Let
$$V_1=e^{ik}\left(
\begin{array}{cc}
\sqrt{p} & 0\\
0 & -\sqrt{p}
\end{array}
\right),
\sp V_2=e^{il}\left(
\begin{array}{cc}
\sqrt{1-p} & 0\\
0 & -\sqrt{1-p}
\end{array}
\right),
$$
where $k,l\in\mathbb{R}$, $p\in (0,1)$. Then $V_1^* V_1+V_2^* V_2=I$. A simple calculation shows that $\rho_2=0$, and then
$$\rho=\left(
\begin{array}{cc}
q & 0\\
0 & 1-q
\end{array}
\right)
$$
is invariant to $\Lambda(\rho)=V_1\rho V_1^*+V_2\rho V_2^*$, for  $q\in (0,1)$.
\end{examp}

\qee

Now we make a few considerations about the Ruelle operator $\mathcal{L}$ defined before. In particular, we show that Perron's classic eigenvalue problem is a particular case of the problem for the operator $\mathcal{L}$ acting on matrices. Let
$$V_1=\left(
\begin{array}{cc}
p_{00} & 0\\
0 & 0
\end{array}
\right),
\sp V_2=\left(
\begin{array}{cc}
0 & p_{01}\\
0 & 0
\end{array}
\right)
$$
$$V_3=\left(
\begin{array}{cc}
0 & 0\\
p_{10} & 0
\end{array}
\right),
\sp V_4=\left(
\begin{array}{cc}
0 & 0\\
0 & p_{11}
\end{array}
\right),
\sp \rho=\left(
\begin{array}{cc}
\rho_1 & \rho_2\\
\rho_3 & \rho_4
\end{array}
\right)
$$
Define
\begin{equation}
\mathcal{L}(\rho)=\sum_{i=1}^4q_i(\rho)V_i\rho V_i^{*}
\end{equation}
We have that $\mathcal{L}(\rho)=\rho$ implies $\rho_2=0$ and
\begin{equation}\label{densist1}
a\rho_1+b\rho_4=\rho_1
\end{equation}
\begin{equation}\label{densist2}
c\rho_1+d\rho_4=\rho_4
\end{equation}
where
$$a=q_1p_{00}^2,\sp b=q_2p_{01}^2,\sp c=q_3p_{10}^2,\sp d=q_4p_{11}^2$$
Solving (\ref{densist1}) and (\ref{densist2}) in terms of $\rho_1$ gives
$$\rho_1=\frac{b}{1-a}\rho_4,\sp \rho_1=\frac{1-d}{c}\rho_4$$
that is,
\begin{equation}\label{condsq}
\frac{b}{1-a}=\frac{1-d}{c}
\end{equation}
which is a restriction over the $q_i$. For simplicity we assume here that the $q_i$ are constant. One can show that
\begin{equation}\label{pirhoc2}
\rho=\left(
\begin{array}{cc}
\frac{q_2p_{01}^2}{q_2p_{01}^2-q_1p_{00}^2+1} & 0\\
0 & \frac{1-q_1p_{00}^2}{q_2p_{01}^2-q_1p_{00}^2+1}
\end{array}
\right)=\left(
\begin{array}{cc}
\frac{1-q_4p_{11}^2}{1-q_4p_{11}^2+q_3p_{10}^2} & 0\\
0 & \frac{q_3p_{10}^2}{1-q_4p_{11}^2+q_3p_{10}^2}
\end{array}
\right)
\end{equation}
Now let
$$P=\sum_i V_i=\left(
\begin{array}{cc}
p_{00} & p_{01}\\
p_{10} & p_{11}
\end{array}
\right),$$ be a column-stochastic matrix. Let $\pi=(\pi_1,\pi_2)$ such that $P\pi=\pi$. Then

\begin{equation}\label{pirhoc1}
\pi=(\frac{p_{01}}{p_{01}-p_{00}+1},\frac{1-p_{00}}{p_{01}-p_{00}+1})\end{equation}

Comparing (\ref{pirhoc1}) and (\ref{pirhoc2}) suggests that we should fix
\begin{equation}\label{solnat}
q_1=\frac{1}{p_{00}},\sp q_2=\frac{1}{p_{01}},\sp q_3=\frac{1}{p_{10}}, \sp q_4=\frac{1}{p_{11}}
\end{equation}
Then the nonzero entries of $\rho$ are equal to the entries of $\pi$ and therefore we associate the fixed point of $P$ to the fixed point of some $\mathcal{L}$ in a natural way. But note that such a choice of $q_i$ is not unique, because
\begin{equation}\label{sol4mat1}
q_2=\frac{1-q_1p_{00}^2}{p_{01}p_{10}},\sp q_4=\frac{1-q_3p_{10}p_{01}}{p_{11}^2},
\end{equation}
for any $q_1,q_3$ also produces $\rho$ with nonzero coordinates equal to the coordinates of $\pi$. We also note that the above calculations can be made by taking the $V_i$ matrices with nonzero entries equal to $\sqrt{p_{ij}}$ instead of $p_{ij}$.

\bigskip

Now we consider the following problem. Let
$$V_1=\left(
\begin{array}{cc}
h_{00} & 0\\
0 & 0
\end{array}
\right),
\sp V_2=\left(
\begin{array}{cc}
0 & h_{01}\\
0 & 0
\end{array}
\right)
,\sp V_3=\left(
\begin{array}{cc}
0 & 0\\
h_{10} & 0
\end{array}
\right)$$
$$V_4=\left(
\begin{array}{cc}
0 & 0\\
0 & h_{11}
\end{array}
\right),
\sp H=\sum_i V_i,
\sp \rho=\left(
\begin{array}{cc}
\rho_1 & \rho_2\\
\rho_3 & \rho_4
\end{array}
\right)
$$

Define
\begin{equation}
\mathcal{L}(\rho)=\sum_{i=1}^4 q_iV_i\rho V_i^*,
\end{equation}
where $q_i\in\mathbb{R}$. Assume that $h_{ij}\in\mathbb{R}$, so we want to obtain $\lambda$ such that
$\mathcal{L}(\rho)=\lambda\rho$, $\lambda\neq 0$, and $\lambda$ is the largest eigenvalue. With a few calculations we obtain $\rho_2=\rho_3=0$,
$$q_1h_{00}^2\rho_1+q_2h_{01}^2\rho_4=\lambda\rho_1$$
$$q_3h_{10}^2\rho_1+q_4h_{11}^2\rho_4=\lambda\rho_4$$
that is,
\begin{equation}\label{yetanother1}
a\rho_1+b\rho_4=\lambda\rho_1
\end{equation}
\begin{equation}\label{yetanother2}
c\rho_1+d\rho_4=\lambda\rho_4,
\end{equation}
with
$$a=q_1h_{00}^2,\sp b=q_2h_{01}^2,\sp c=q_3h_{10}^2,\sp d=q_4h_{11}^2$$
Therefore
$$\rho=\left(
\begin{array}{cc}
\frac{\lambda-d}{c}\rho_4 & 0\\
0 & \rho_4
\end{array}
\right)=\left(
\begin{array}{cc}
\frac{b}{\lambda-a}\rho_4 & 0\\
0 & \rho_4
\end{array}
\right)
$$
and
$$\frac{\lambda-d}{c}=\frac{b}{\lambda-a}$$
Solving for $\lambda$, we obtain the eigenvalues
$$\lambda=\frac{a+d}{2}\pm\frac{\zeta}{2}=\frac{a+d}{2}\pm\frac{\sqrt{(d-a)^2+4bc}}{2}$$
$$=\frac{1}{2}\Big(q_1h_{00}^2+q_4h_{11}^2\pm\sqrt{(q_4h_{11}^2-q_1h_{00}^2)^2+4q_2q_3h_{01}^2h_{10}^2}\Big),$$
where
$$\zeta=\sqrt{(d-a)^2+4bc}=\sqrt{(q_4h_{11}^2-q_1h_{00}^2)^2+4q_2q_3h_{01}^2h_{10}^2}$$
and the associated eigenfunctions
$$\rho=\left(
\begin{array}{cc}
\frac{a-d\pm\zeta}{2c}\rho_4 & 0\\
0 & \rho_4
\end{array}
\right)=\left(
\begin{array}{cc}
\frac{2b}{d-a\pm\zeta}\rho_4 & 0\\
0 & \rho_4
\end{array}
\right)
$$
But $\rho_1+\rho_4=1$ so we obtain
$$
\rho=\left(
\begin{array}{cc}
\frac{a-d\pm\zeta}{a-d\pm\zeta+2c} & 0\\
0 & \frac{2c}{a-d\pm\zeta+2c}
\end{array}
\right)$$
\begin{equation}\label{eqcomplvejaa}
=\left(
\begin{array}{cc}
\frac{q_1h_{00}^2-q_4h_{11}^2\pm\zeta}{q_1h_{00}^2-q_4h_{11}^2\pm\zeta+2q_3h_{10}^2} & 0\\
0 & \frac{2q_3h_{10}^2}{q_1h_{00}^2-q_4h_{11}^2\pm\zeta+2q_3h_{10}^2}
\end{array}
\right)
\end{equation}
that is,
$$
\rho=\left(
\begin{array}{cc}
\frac{-2b}{a-2b-d\mp\zeta} & 0\\
0 & \frac{a-d\mp\zeta}{a-2b-d\mp\zeta}
\end{array}
\right)$$
\begin{equation}\label{eqcomplveja}
=\left(
\begin{array}{cc}
\frac{-2q_2h_{01}^2}{q_1h_{00}^2-2q_2h_{01}^2-q_4h_{11}^2\mp\zeta} & 0\\
0 & \frac{q_1h_{00}^2-q_4h_{11}^2\mp\zeta}{q_1h_{00}^2-2q_2h_{01}^2-q_4h_{11}^2\mp\zeta}
\end{array}
\right)
\end{equation}
Therefore we obtained that $\rho_1,\rho_4, q_1,\dots,q_4, \lambda$ are implicit solutions for the set of equations (\ref{yetanother1})-(\ref{yetanother2}). Recall that in this case we obtained $\rho_2=\rho_3=0.$

\bigskip
Now we consider the problem of finding the eigenvector associated to the dominant eigenvalue of $H$. The eigenvalues are
$$\lambda=\frac{1}{2}\Big(h_{00}+h_{11}\pm\sqrt{(h_{00}-h_{11})^2+4h_{01}h_{10}}\Big)$$
Then we can find $v$ such that $Hv = \lambda v$ from the set of equations
\begin{equation}\label{yetanother3}
h_{00}v_1+h_{01}v_2=\lambda v_1
\end{equation}
\begin{equation}\label{yetanother4}
h_{10}v_1+h_{11}v_2=\lambda v_2
\end{equation}
which determine $v_1, v_2,\lambda$ implicitly. Note that if we set

\begin{equation}
q_1=\frac{1}{p_{00}},\sp q_2=\frac{1}{p_{01}},\sp q_3=\frac{1}{p_{10}}, \sp q_4=\frac{1}{p_{11}}
\end{equation}
we have that the set of equations (\ref{yetanother1})-(\ref{yetanother2}) and
(\ref{yetanother3})-(\ref{yetanother4}) are the same.
Hence we conclude that Perron's classic eigenvalue
problem is a particular case of the problem for
$\mathcal{L}$ acting on matrices.

\qee

A different analysis in the quantum setting which is related to Perron's
theorem is presented in \cite{bruzda}.

\section{A theorem on eigenvalues for the Ruelle operator}\label{sec_eigenvalues}

The following proposition is inspired in \cite{par}. We say that a hermitian operator  $P:V\to V$ on a Hilbert space $(V,\langle\cdot\rangle)$ is {\bf positive} if $\langle Pv,v\rangle\geq 0$, for all $v\in V$, denoted $P\geq 0$. Consider the positive operator $\mathcal{L}_{W,V}:\mathcal{PH}_N\to \mathcal{PH}_N$,
\begin{equation}
\mathcal{L}_{W,V}(\rho):=\sum_{i=1}^k tr(W_i\rho W_i^*)V_i\rho V_i^*.
\end{equation}

We point out that this operator is completely general. In an analogy with the classical case we can say it corresponds to the general Perron Theorem for positive matrices (having positive eigenvalues which can be bigger or smaller than one), by the other hand  the setting described in
\cite{lozinski}, \cite{wsbook} "basically" considers the analogous case of the Perron Theorem for stochastic matrices.

We need a result in this form in order to better understand the Pressure problem which will be described later.

\begin{pro}\label{proavl1}
There exists $\rho\in\mathcal{M}_N$ and $\beta>0$ such that $\mathcal{L}_{W,V}(\rho)=\beta\rho$. The value $\beta$ is obtained explicitly: $\beta=tr(\mathcal{L}_{W,V}(\rho))$.
\end{pro}
{\bf Proof} Define $\mathcal{L}_n:\mathcal{M}_N\to\mathcal{M}_N$,
$$\mathcal{L}_n(\rho):=\frac{\mathcal{L}_{W,V}(\rho+\frac{I}{n})}{tr(\mathcal{L}_{W,V}(\rho+\frac{I}{n}))} \sp , \sp n\geq 1$$
The operator above is well defined. In fact, note that $\mathcal{L}_{W,V}(\rho)$, $W_jW_j^*$, $V_jV_j^*$ are positive for all $j$. Then
$$tr\Big[\sum_i tr\Big(W_i(\rho+\frac{I}{n})W_i^*\Big)V_i(\rho+\frac{I}{n})V_i^*\Big]
=\sum_i tr\Big(W_i(\rho+\frac{I}{n})W_i^*\Big)tr(V_i(\rho+\frac{I}{n})V_i^*)$$
$$=\sum_i tr(W_i\rho W_i^*+\frac{1}{n}W_iW_i^*)tr(V_i\rho V_i^*+\frac{1}{n}V_iV_i^*)\geq$$
$$\geq \sum_i tr(W_i\rho W_i^*)tr(V_i\rho V_i^*)=tr(\mathcal{L}_{W,V})$$
We know that for any positive operator $P\neq 0$, if $\{v_1,\dots, v_N\}$ is a orthonormal base for $\mathcal{H}_N$, then
$$tr(P)=\sum_{i=1}^N \langle Pv_i,v_i\rangle>0$$
Therefore, $tr(\mathcal{L}_{W,V}(\rho+\frac{I}{n}))>0$, $n\geq 1$. Hence $\mathcal{L}_n(\rho)$ is well defined.

\bigskip

We know that $\mathcal{M}_N$ is compact and convex, so we can apply Schauder's theorem for each of the mappings $\mathcal{L}_n$, $n\geq 1$ and get $\rho_n\in\mathcal{M}_N$ such that $$\mathcal{L}_n(\rho_n)=\rho_n\sp \Rightarrow\sp\mathcal{L}_{W,V}(\rho_n+\frac{I}{n})=\beta_n\rho_n,\sp n\geq 1$$
where
$$\beta_n:=tr(\mathcal{L}_{W,V}(\rho_n+\frac{I}{n}))$$
By the compacity of $\mathcal{M}_N$, we can choose a point $\rho\in\mathcal{M}_N$ which is limit of the sequence $\{\rho_n\}_{n=1}^\infty$ and then, by continuity, $\mathcal{L}_{W,V}(\rho)=\beta\rho$, where $\beta=tr(\mathcal{L}_{W,V}(\rho))$. Also, note that $\beta \geq 0$, because if $\{v_1,\dots, v_N\}$ is a orthonormal base of $\mathcal{H}_N$,
$$tr(\mathcal{L}_{W,V}(\rho))=\sum_{i=1}^N \langle\mathcal{L}_{W,V}(\rho) v_i,v_i\rangle\geq 0,$$
since $\mathcal{L}_{W,V}(\rho)$ is positive, and the inequality will be equal to zero if and only if  $\mathcal{L}_{W,V}(\rho)$ is the zero operator. Hence, we proved that there exists $\rho\in\mathcal{M}_N$ and $\beta>0$ such that $\mathcal{L}_{W,V}(\rho)=\beta\rho$.

\qed

\section{Vector integrals and barycenters}\label{v_bari}

We recall here a few basic definitions. For more details, see \cite{lozinski} and \cite{wsbook}. Let $X$ be a metric space. Let $(V,+,\cdot)$ be a real vector space, and $\tau$ a topology on $V$. We say that $(V,+,\cdot;\tau)$ is a topological vector space if it is Hausdorff and if the operations $+$ and $\cdot$ are continuous. For instance, in the context of density matrices, we will consider $V$ as the space of hermitian operators $\mathcal{H}_N$ and $X$ will be the space of density matrices $\mathcal{M}_N$.

\begin{defi}
    Let $(X,\Sigma)$ be a measurable space, let $\mu\in M(X)$, let  $(V,+,\cdot ; \tau)$ be a locally convex space and let $f:X\to V$. we say that $x\in V$ is the {\bf integral} of $f$ in $X$, denoted by \begin{equation}
    x:=\int_X fd\mu
    \end{equation}
if
\begin{equation}
\Psi(x)=\int_X\Psi\circ fd\mu,
\end{equation}
for all $\Psi\in V^*$.
\end{defi}

It is known that if we have a compact metric space $X$, $V$ is a
locally convex space and $f:X\to V$ is a continuous function such
that $\ov{co}f(X)$ is compact then the integral of $f$ in $X$
exists and belongs to $\ov{co}f(X)$. We will also use the
following well-known result, the barycentric formula:
\begin{pro}\label{winkler_teo}
\cite{winkler} Let $V$ be a locally convex space, let $E\subset V$ be a complete, convex and bounded set, and $\mu\in M^1(E)$.
Then there is a unique $x\in E$ such that
$$l(x)=\int_E l d\mu,$$
for all $l\in V^*$.
\end{pro}
In the context of QIFS, we can take $V=E=\mathcal{M}_N$.

\section{Example: density matrices}\label{sec_densities}

In this section we briefly review how the constructions of the previous section adjust to the case of density matrices. Define $V:=\mathcal{H}_N$, $V^+:=\mathcal{PH}_N$ (note that such space is a convex cone), and let the partial order $\leq$ on $\mathcal{PH}_N$ be $\rho\leq\psi$ if and only if $\psi-\rho\geq 0$, i.e., if $\psi-\rho$ is positive. Then
$$(V,V^+,e)=(\mathcal{H}_N,\mathcal{PH}_N,tr),$$
is a regular state space \cite{wsbook}. Also, the set $B$ of unity trace in $V^+$ is, of course, the space of density matrices, so $B=\mathcal{M}_N$.

\bigskip
Let $Z\subset V^*$ be a nonempty vector subspace of $V^*$. The smallest topology in $V$ such that every functional defined in $Z$ is continuous on that topology, denoted by $\sigma(V,Z)$, turns $V$ into a locally convex space. In particular, $\sigma(V,V^*)$ is the weak topology in $V$. If $(V,\Vert\cdot\Vert)$ is a normed space, then $\sigma(V^*,V)$ is called a weak$^*$ topology in $V^*$ (we identify $V$ with a subspace of $V^{**})$. We also have that $(C,\tau)=(\mathcal{PH}_N,\tau)$, where $\tau$ is the weak$^*$ topology (and which is equal to the Euclidean, see \cite{wsbook}) is a metrizable compact structure. In this case we have that $B_C=B\cap C=\mathcal{M}_N$.

\begin{defi}
A {\bf Markov operator} for probability measures is an operator $P:M^1(X)\to M^1(X)$ such that
\begin{equation}
P(\lambda\mu_1+(1-\lambda)\mu_2)=\lambda P\mu_1+(1-\lambda)P\mu_2,
\end{equation}
for $\mu_1,\mu_2\in M^1(X)$, $\lambda\in (0,1)$.
\end{defi}
An example of such operator is the one which we have defined before and we denote it by $\mathcal{V} : M^1(\mathcal{M}_N)\to M^1(\mathcal{M}_N)$,
\begin{equation}\label{markov_induced}
(\mathcal{V}\nu)(B)=\sum_{i=1}^k\int_{F_i^{-1}(B)}p_id\nu
\end{equation}
We call it the Markov operator induced by the QIFS $\mathcal{F}=\{\mathcal{M}_N,F_i,p_i\}_{i=1,\dots, k}$. Define
$$m_b(X):=\{f:X\to\mathbb{R} : \textrm{f is bounded, measurable}\}$$
Then define $\mathcal{U}:m_b(X)\to m_b(X)$,
\begin{equation}\label{v_dual_op}
(\mathcal{U}f)(x):=\sum_{i=1}^k p_i(x)f(F_i(x))
\end{equation}

\begin{pro}\label{pdualbas}\cite{wsbook}
Let $f\in m_b(X)$ and $\mu\in M^1(X)$, then
\begin{equation}
\langle f,\mathcal{V}\mu\rangle=\langle\mathcal{U}f,\mu\rangle=\sum_{i=1}^k\int p_i(f\circ F_i)d\mu,
\end{equation}
where $\langle f,\mu\rangle$ denotes the integral of $f$ with respect to $\mu$.
\end{pro}

\begin{defi} An operator $Q:V^+\to V^+$ is {\bf submarkovian} if
\begin{enumerate}
\item $Q(x+y)=Q(x)+Q(y)$
\item $Q(\alpha x)=\alpha Q(x)$
\item $\Vert Q(x)\Vert\leq\Vert x\Vert,$
\end{enumerate}
for all $x$, $y\in V^+$, $\alpha >0$.
\end{defi}
Every submarkovian operator $Q:V^+\to V^+$ can be extended in a unique way to a positive linear contraction on $V$, see \cite{wsbook}.

\begin{defi} Let $P:V^+\to V^+$ a Markov operator and let
$P_i:V^+\to V^+$, $i=1,\dots, k$ be submarkovian operators such that $P=\sum_i P_i$. We say that $(P,\{P_i\}_{i=1}^k)$ is a {\bf Markov pair}.
\end{defi}
From \cite{wsbook}, we know that there is a 1-1 correspondence between homogeneous IFS and Markov pairs.

\section{Some lemmas for IFS}\label{sec_somelemmas}

We want to understand the structure of $\Lambda:\mathcal{M}_N\to \mathcal{M}_N$,
\begin{equation}
\Lambda(\rho):=\sum_{i=1}^k p_iF_i=\sum_{i=1}^k tr(W_i\rho W_i^*)\frac{V_i\rho V_i^*}{tr(V_i\rho V_i^*)},
\end{equation}
where $V_i$, $W_i$ are linear, $\sum_i W_i^*W_i=I$. Such operator is associated in a natural way to an IFS which is not homogeneous. In this section we state a few useful properties which are relevant for our study. The following lemmas hold for any IFS, except for lemma \ref{lema_imp01}, where the proof presented here is valid only for homogeneous IFS.

\begin{lem}\label{umlemautil1}
Let $\{X,F_i,p_i\}_{i=1,\dots, k}$ be an IFS, $\Psi$ a linear functional on $X$. Then $\mathcal{U}\circ\Psi=\Psi\circ\Lambda$,  where $\mathcal{U}$ is given by $(\ref{v_dual_op})$.
\end{lem}
{\bf Proof} We have
$$(\mathcal{U}\Psi)(x)=\sum_{i}p_i(x)\Psi(F_i(x))=\Psi(\sum_i p_i(x)F_i(x))=\Psi(\Lambda(x))$$

\qed

\begin{cor}\label{ccc1}
Let $\mathcal{F}=(X,F_i,p_i)_{i=1,\dots, k}$ be an IFS and let $\rho_0\in X$. Then $\Lambda(\rho_0)=\rho_0$ if and only if $\mathcal{U}(\Psi(\rho_0))=\Psi(\rho_0)$, for all $\Psi$ linear functional.
\end{cor}
{\bf Proof} Suppose that $\mathcal{L}(\rho_0)=\rho_0$. Then
$$\mathcal{U}(\Psi(\rho_0))=\sum_i p_i(\rho_0)\Psi(F_i(\rho_0))=\Psi(\sum_i p_i(\rho_0)F_i(\rho_0))=\Psi(\Lambda(\rho_0))=\Psi(\rho_0)$$
Conversely, if $\mathcal{U}(\Psi(\rho_0))=\Psi(\rho_0)$, then
$$\Psi(\Lambda(\rho_0))=\mathcal{U}(\Psi(\rho_0))=\Psi(\rho_0)$$

\qed

\begin{lem}\label{lemabas1}
Let $\mathcal{F}=\{X, F_i, p_i\}_{i=1,\dots, k}$ be an IFS.
\begin{enumerate}
\item Let $\rho_0\in X$ such that $F_i(\rho_0)=\rho_0$, $i=1,\dots, k$. Then $\mathcal{V}\delta_{\rho_0}=\delta_{\rho_0}$.
\item Let $\rho_0\in X$ such that $\mathcal{V}\delta_{\rho_0}=\delta_{\rho_0}$, then $\Lambda(\rho_0)=\rho_0$.
\end{enumerate}
\end{lem}
{\bf Proof} 1. We have
$$\mathcal{V}\delta_{\rho_0}(B)=\sum_{i=1}^k\int_{F_i^{-1}(B)}p_id\delta_{\rho_0}=\sum_{i=1}^k\int p_i(\rho)1_B(F_i(\rho))d\delta_{\rho_0}$$
$$=\sum_{i=1}^kp_i(\rho_0)1_B(F_i(\rho_0))=\sum_{i=1}^kp_i(\rho_0)1_B(\rho_0)=\delta_{\rho_0}(B)$$
2. Let $\Psi$ be a linear functional. Then
$$\Psi(\Lambda(\rho_0))=\mathcal{U}(\Psi(\rho_0))=\int\mathcal{U}(\Psi(\rho))d\delta_{\rho_0}=\int \Psi(\rho)d\mathcal{V}\delta_{\rho_0}$$
$$=\int \Psi(\rho)d\delta_{\rho_0}=\Psi(\rho_0)$$

\qed

\begin{lem}\label{lema_imp01}
Let $\{X,F_i,p_i\}_{i=1,\dots, k}$ be a homogeneous IFS, $\Lambda=\sum_i p_i F_i$.

\begin{enumerate}
\item Let $\rho_{\nu}$ be the barycenter of a probability measure $\nu$. Then $\Lambda(\rho_{\nu})$ is the barycenter of $\mathcal{V}\nu$, where $\mathcal{V}$ is the associated Markov operator.

\item Let $\mu$ be an invariant probability measure for $\mathcal{V}$. Then the barycenter of $\mu$, denoted by $\rho_\mu$, is a fixed point of $\Lambda$.

\end{enumerate}

\end{lem}
{\bf Proof} 1. We have, for $\Psi$ linear functional,
$$\Psi(\Lambda(\rho_{\nu}))=\int\Psi(\Lambda(\rho))d\nu=\int \mathcal{U}\circ \Psi d\nu=\int\Psi d\mathcal{V}\nu$$

\bigskip

2. By lemma (\ref{umlemautil1}), we have
$$\Psi(\Lambda(\rho_{\mu}))=\mathcal{U}\circ\Psi(\rho_{\mu})=\int\mathcal{U}\circ\Psi d\mu=\int\Psi d\mathcal{V}\mu=\int\Psi d\mu=\Psi(\rho_{\mu}),$$
where the fact that $\mathcal{U}\circ\Psi$ is linear follows from the homogeneity of $\mathcal{F}$.

\qed

In order to prove uniqueness in  item (2) above it would be necessary to assume hyperbolicity \cite{slom}. It is known that without this hypothesis even in the classical case (for transformations for instance) it can happen the phenomena of phase transition (two or more probabilities which are solutions) \cite{young} \cite{FL}. The present setting contains the classical case and therefore in general there is no uniqueness.

\begin{examp}\label{exemplo_ifs} Let $k=N=2$,
$$V_1=\left(
\begin{array}{cc}
-1 & 0 \\
0 & 1
\end{array}
\right)
,\sp
V_2=\left(
\begin{array}{cc}
0 & -\frac{3\sqrt{2}}{4} \\
-\frac{3\sqrt{2}}{2} & 0
\end{array}
\right),$$
$W_1=(1/2)I$, $W_2=(\sqrt{3}/2)I$. Then
$$\Lambda(\rho)=\sum_i p_i(\rho)F_i(\rho)=\sum_i tr(W_i\rho W_i^*)\frac{V_i\rho V_i^*}{tr(V_i\rho V_i^*)}$$
$$=\frac{1}{4}V_1\rho V_1^*+\frac{3}{4}\frac{V_2\rho V_2^*}{tr(V_2\rho V_2^*)}
=\frac{1}{4}V_1\rho V_1^*+\frac{3}{4}\frac{V_2\rho V_2^*}{(\frac{9}{8}+\frac{27}{8}\rho_1)}$$
induces an IFS and it is such that  $\rho_0=\frac{1}{3}\vert 0\rangle\langle 0\vert+\frac{2}{3}\vert 1\rangle\langle 1\vert$ is a fixed point, with $F_1(\rho_0)=F_2(\rho_0)=\rho_0$. We can apply lemma \ref{lemabas1} and conclude that $\delta_{\rho_0}$ is an invariant measure for the Markov operator $\mathcal{V}$ associated to the IFS determined by $p_i$ and $F_i$.
\end{examp}

\qee

The following lemma, a simple variation from results seen in
\cite{wsbook}, specifies a condition we need in
order to obtain a fixed point for $\Lambda$ from a certain
measure which is invariant for the Markov operator $\mathcal{V}$.

\begin{lem}\label{ptofixapa}
Let $\{\mathcal{M}_N,F_i,p_i\}_{i=1,\dots, k}$ be an IFS which admits an attractive invariant measure $\mu$ for $\mathcal{V}$. Then $\lim_{n\to\infty}\Lambda^n(\rho_0)=\rho_{\mu}$, for every $\rho_0\in\mathcal{M}_N$, where $\rho_{\mu}$ is the barycenter of $\mu$.
\end{lem}

{\bf Proof} Let $\rho_0\in\mathcal{M}_N$. Then
$$\Psi(\Lambda^n(\rho_0))=\mathcal{U}^n(\Psi(\rho_0))=\int \mathcal{U}^n(\Psi(\rho))d\delta_{\rho_0}=\int \Psi(\rho)d\mathcal{V}^n\delta_{\rho_0}$$
so $\Psi(\Lambda^n(\rho_0))\to\int\Psi(\rho)d\mu=\Psi(\rho_{\mu})$, as $n\to\infty$, for all $\Psi$ linear functional. Hence, $\Lambda^n(\rho_0)\to\rho_{\mu}$ as $n\to\infty$, for all $\rho_0\in\mathcal{M}_N$.
\qed

In lemma \ref{ptofixapa}, we have a general QIFS and an attractive invariant $\mu$, then $\mu$ is the unique invariant measure, an easy consequence of  attractivity \cite{wsbook}. In general, we will be interested in QIFS which has an attractive invariant measure. This will follow if we assume hyperbolicity.

\section{Integral formulae for the entropy of IFS}\label{sec_entr_ifs}

Part of the results we present here in this section are variations of results presented in \cite{wsbook}. Let $(X,d)$ be a complete separable metric space. Let $(V,V^+,e)$ be a complete state space, $B=\{x\in V^+: e(x)=1\}$ and $\mathcal{F}=(X,F_i,p_i)_{i=1,\dots, k}$ the homogeneous IFS induced by the Markov pair $(\Lambda,\{\Lambda_i\}_{i=1}^k)$. Now define
$I_{k}:=\{1,\dots, k\}$. Let $n\in\mathbb{N}$, $\iota\in I_k^n$, $i\in I_k$. Define
$F_{\iota i}:=F_i\circ F_\iota$ and
\begin{equation}
p_{\iota i}(x)=\left\{\begin{array}{ll}
p_i(F_{\iota}x)p_{\iota}(x) & \textrm{ if } p_\iota(x)\neq 0\\
0 & \textrm{ otherwise }
\end{array}\right.
\end{equation}

\begin{pro}\label{un_etc}
Let $n\in\mathbb{N}$, $f\in m_b(X)$, $x\in X$. Then
$$(\mathcal{U}^n f)(x)=\sum_{\iota\in I_k^n} p_\iota(x)f(F_\iota(x))$$
\end{pro}

\begin{pro}
Let $x\in B$, $n\in\mathbb{N}$. Then
$$\Lambda^n(x)=\sum_{\iota\in I_k^n} p_\iota(x)F_\iota(x).$$
\end{pro}

\begin{pro}\label{ws_ident}
Let $\mathcal{F}$ be an IFS and let $g:B\to\mathbb{R}$. Then for $n\in\mathbb{N}$,
\begin{enumerate}
\item If g is concave (resp. convex, affine) then $\mathcal{U}^n g\leq g\circ \Lambda^n$ (resp. $\mathcal{U}^n g\geq g\circ \Lambda^n$, $\mathcal{U}^n g=g\circ\Lambda^n$).
\item If $\ov{x}$ is a fixed point for $\Lambda$ then the sequence $(\mathcal{U}^ng)(\ov{x}))_{n\in\mathbb{N}}$ is decreasing (resp. increasing, constant) if $g$ is concave (resp. convex, affine).

\bigskip
Also suppose that $\mathcal{F}$ is homogeneous. Then

\item If g is concave (resp. convex, affine), then $\mathcal{U}g$ is concave (resp. convex, affine).
\end{enumerate}
\end{pro}

\bigskip
We recall some well-known definitions and results. Define
$\eta:\mathbb{R}^+\to\mathbb{R}$ as
\begin{equation}
\eta(x)=\left\{\begin{array}{ll}
-x\log{x} & \textrm{ if } x\neq 0\\
0 & \textrm{ if } x=0
\end{array}\right.
\end{equation}
Then the Shannon-Boltzmann entropy function is $h:X\to \mathbb{R}^+$,
\begin{equation}
h(x):=\sum_{i=1}^k\eta(p_i(x))
\end{equation}
Let $n\in\mathbb{N}$. Define the {\bf partial entropy} $H_n:X\to\mathbb{R}^+$ as
\begin{equation}
H_n(x):=\sum_{\iota\in I_k^n}\eta(p_\iota (x)),
\end{equation}
for $n\geq 1$ and $H_0(x):=0$, $x\in X$. Define, for $x\in X$,
\begin{equation}
\ov{\mathcal{H}}(x):=\limsup_{n\to\infty}\frac{1}{n}H_n(x),
,\sp\underline{\mathcal{H}}(x):=\liminf_{n\to\infty}\frac{1}{n}H_n(x),
\end{equation}
the {\bf upper} and {\bf lower} entropy on x. If such limits are equal, we call its common value the {\bf entropy on x}, denoted by $\mathcal{H}(x)$.

\bigskip
Denote by $M^{\mathcal{V}}(X)$ the set of $\mathcal{V}$-invariant
probability measures on $X$. Let $\mu\in M^{\mathcal{V}}(X)$. The
{\bf partial entropy of the measure $\mu$} is defined by
\begin{equation}
H_n(\mu):=\sum_{\iota\in I_k^n}\eta(\langle p_\iota,\mu\rangle),
\end{equation}
for $n\geq 1$ and $H_0(\mu):=0$.

\begin{pro}
Let $\mu\in M^{\mathcal{V}}(X)$. Then the sequences
$(\frac{1}{n}H_n(\mu))_{n\in\mathbb{N}}$ and
$(H_{n+1}(\mu)-H_{n}(\mu))_{n\in\mathbb{N}}$ are nonnegative,
decreasing, and have the same limit.
\end{pro}
We denote the common limit of the sequences mentioned in the proposition above
as $\mathcal{H}(\mu)$ and we call it the {\bf entropy of the measure} $\mu$, i.e.,
\begin{equation}
\mathcal{H}(\mu):=\lim_{n\to\infty}\frac{1}{n}H_n(\mu)=\lim_{n\to\infty}(H_{n+1}(\mu)-H_{n}(\mu))
\end{equation}

\bigskip

The following result gives us an integral formula for entropy, and
also a relation between the entropies defined before. We write
$S(\mu):=M^{\mathcal{V}}(X)\cap
\textrm{Lim}(\mathcal{V}^n\mu)_{n\in\mathbb{N}},$ where
$\textrm{Lim}(\mathcal{V}^n\mu)_{n\in\mathbb{N}}$ is the convex
hull of the set of accumulation points of
$(\mathcal{V}^n\mu)_{n\in\mathbb{N}}$, and $S_{\mathcal{F}}(\mu)$
is the set $S(\mu)$ associated to the Markov operator induced by
the IFS $\mathcal{F}$. For the definition of compact structure and
$(C,\tau)$-continuity, see \cite{wsbook}.

\begin{teo}\label{teo71entropia}
\cite{wsbook} (Integral formula for entropy of homogeneous IFS, compact case). Let $(C,\tau)$ be a metrizable compact structure $(V,V^+,e)$ such that $(\Lambda,\{\Lambda_i\}_{i=1}^k)$ is $(C,\tau)$-continuous. Assume that  $\rho_0\in B_C:=B\cap C$ is such that $\Lambda(\rho_0)=\rho_0$. Then
$$\mathcal{H}(\rho_0)=\mathcal{H}(\nu)=\int_X hd\nu$$
for each $\nu\in S_{\mathcal{F}_C}(\delta_{\rho_0})$, where $\mathcal{F}_C$ is the IFS $\mathcal{F}$ restricted to $(B_C,\tau)$.
\end{teo}

The analogous result for hyperbolic IFS is the following.

\begin{teo}\label{teo71entropiab}
\cite{wsbook} Let $\mathcal{F}=(X,F_i,p_i)_{i=1,\dots, k}$ be a hyperbolic IFS, $x\in X$, $\mu\in M^1(X)$ an attractive invariant measure for $\mathcal{F}$. Then
$$\mathcal{H}(x)=\lim_{n\to\infty}(H_{n+1}(x)-H_n(x))$$
and
$$\mathcal{H}(x)=\mathcal{H}(\mu)=\int_X hd\mu .$$
\end{teo}

\section{Some calculations on entropy}\label{sec_calc_entr}

Let $U$ be a unitary matrix of order $mn$ acting on $\mathcal{H}_m\otimes\mathcal{H}_n$. Its Schmidt decomposition is
$$U=\sum_{i=1}^K\sqrt{q_i}V_i^A\otimes V_i^B,\sp K=min\{m^2,n^2\}$$ The operators $V_i^A$ and $V_i^B$ act on certain Hilbert spaces $\mathcal{H}_m$ and $\mathcal{H}_n$, respectively. We also have that  $\sum_{i=1}^K q_i=1$. Let $\sigma=\rho_A\otimes\rho_*^B=\rho_A\otimes I_n/n$ and define
$$\Lambda(\rho_A):=tr_B(U\sigma U^*)=\sum_{i=1}^Kq_i V_i^A\rho_A V_i^{A*}$$
Above, recall that the partial trace is
$$tr_B(\vert a_1\rangle\langle a_2\vert\otimes \vert b_1\rangle\langle b_2\vert):
=\vert a_1\rangle\langle a_2\vert tr(\vert b_1\rangle\langle b_2\vert)$$
where $\vert a_1\rangle$ and $\vert a_2\rangle$ are vectors on the state space of $A$ and $\vert b_1\rangle$ and $\vert b_2\rangle$ are vectors on the state space of $B$. The trace on the right side is the usual trace on $B$. A calculation shows that if $\rho_*^A=I_m/m$, then $\Lambda(\rho_*^A)=\rho_*^A$ and so $\Lambda$ is such that $\Lambda(I_m/m)=I_m/m$ and $\Lambda$ is trace preserving.

\bigskip

Let $\mathcal{F}$ be the homogeneous IFS associated to the $V_i^A$, that is,  $p_i(\rho)=tr(q_iV_i^A\rho V_i^{A*})$, $F_i(\rho)=(q_iV_i^A\rho V_i^{A*})/tr(q_iV_i^A\rho V_i^{A*})$ and let $\rho_0$ be a fixed point of $\Lambda=\sum_i p_iF_i$. Following \cite{wsbook}, we have that $\rho_0$ is the barycenter of $\mathcal{V}^n\delta_{\rho_0}$, $n\in\mathbb{N}$. By theorem \ref{teo71entropia}, we can calculate the entropy of such IFS. In this case we have
\begin{equation}\label{expressao_entr}
\mathcal{H}(\rho_0)=\mathcal{H}(\nu)=\int_{\mathcal{M}_N} h d\nu,
\end{equation}
where $\nu\in M^{\mathcal{V}}(X)\cap
\textrm{Lim}(\mathcal{V}^n\delta_{\rho_0})_{n\in\mathbb{N}}$.

\qee

Let $\mathcal{F}=(\mathcal{M}_N,F_i,p_i)_{i=1,\dots, k}$ be an IFS, $\Lambda(\rho)=\sum_i p_iF_i$. Let $\mathcal{U}$ be the conjugate of $\mathcal{V}$. By proposition \ref{un_etc},
$$(\mathcal{U}^n h)(\rho)=\sum_{\iota\in I_k^n(\rho)} p_\iota(\rho)h(F_\iota(\rho))$$
and since $h(\rho)=\sum_{j=1}^k\eta(p_j(\rho))$, we have, for $\iota=(i_1,\dots, i_n)$, and every $\rho_0\in\mathcal{M}_N$,
\begin{equation}\label{inicio_copia}
\int_{\mathcal{M}_N} h d\mathcal{V}^n\delta_{\rho_0}=\int_{\mathcal{M}_N} \mathcal{U}^n h d\delta_{\rho_0}
\end{equation}
\begin{equation}
=-\int_{\mathcal{M}_N}\sum_{\iota\in I_k^n}p_\iota(\rho)\sum_{j=1}^k p_j(F_\iota(\rho))\log{p_j(F_\iota(\rho))}d\delta_{\rho_0}
\end{equation}
\begin{equation}
=-\sum_{\iota\in I_k^n}p_\iota(\rho_0)\sum_{j=1}^k p_j(F_\iota(\rho_0))\log{p_j(F_\iota(\rho_0))}
\end{equation}
\begin{equation}
=-\sum_{\iota\in I_k^n}p_{i_1}(\rho_0)p_{i_2}(F_{i_1}\rho_0)\cdots p_{i_n}(F_{i_{n-1}}(F_{i_{n-2}}(\cdots (F_{i_1}\rho_0))))\times
\end{equation}
\begin{equation}\label{ac_conta1}
\times\sum_{j=1}^k p_j(F_{i_n}(F_{i_{n-1}}(\cdots(F_{i_1}\rho_0))))\log{p_j(F_{i_n}(F_{i_{n-1}}(\cdots(F_{i_1}\rho_0))))}=(\mathcal{U}^n h)(\rho_0)
\end{equation}

Suppose $\Lambda(\rho_0)=\rho_0$. We have by proposition \ref{ws_ident}, since $h$ is concave, that $(\mathcal{U}^nh)_{n\in\mathbb{N}}$ is decreasing, $\mathcal{U}^nh\leq h\circ\Lambda^n$ and so
\begin{equation}\label{ac_conta2}
\int_{\mathcal{M}_N} h d\mathcal{V}^n\delta_{\rho_0}\leq h(\Lambda^n(\rho_0))= h(\rho_0),
\end{equation}
for every $n$.

\section{An expression for a stationary entropy}\label{sec_novaentr}

In this section we present a definition of entropy which captures a stationary behavior. Let $H$ be a hermitian operator and $V_i$, $i=1,\dots, k$ linear operators. We can define the dynamics $F_i:\mathcal{M}_N\to \mathcal{M}_N$:
\begin{equation}\label{din01}
F_i(\rho):=\frac{V_i\rho V_i^*}{tr(V_i\rho V_i^*)}
\end{equation}
Let $W_i$, $i=1,\dots, k$ be linear and such that $\sum_{i=1}^k W_i^*W_i=I$. This determines functions $p_i:\mathcal{M}_N\to\mathbb{R}$,
\begin{equation}\label{prob01}
p_i(\rho):=tr(W_i\rho W_i^*)
\end{equation}
Then we have $\sum_{i=1}^k p_i(\rho)=1$, for every $\rho$. Therefore a family $W:=\{W_i\}_{i=1,\dots, k}$ determines a QIFS $\mathcal{F}_{W}=\{\mathcal{M}_N,F_i,p_i\}_{i=1,\dots, k}$,
with $F_i$, $p_i$ given by (\ref{din01}) and (\ref{prob01}). We introduce the following definition.

\begin{defi} Let $\mathcal{F}_W$ be a QIFS such that there is a unique attractive invariant measure for the associated Markov operator $\mathcal{V}$. Let $\rho_W$ be the barycenter of such measure. Define the {\bf QIFS entropy}:
\begin{equation}\label{nossa_entropia}
h_V(W):=-\sum_{i=1}^k tr(W_i\rho_W W_i^*)\sum_{j=1}^k tr\Big(\frac{W_j V_i\rho_W V_i^* W_j^*}{tr(V_i\rho_W V_i^*)}\Big)\log{tr\Big(\frac{W_j V_i\rho_W V_i^* W_j^*}{tr(V_i\rho_W V_i^*)}\Big)}
\end{equation}

\end{defi}

\bigskip

Remember that by lemma \ref{ptofixapa}, we have that $\rho_W$ is a fixed point for
\begin{equation}\label{lhat1}
\Lambda(\rho)=\Lambda_{\mathcal{F}_W}(\rho):=\sum_{i=1}^k p_i(\rho)F_i(\rho)=\sum_{i=1}^k tr(W_i\rho W_i^*)\frac{V_i\rho V_i^*}{tr(V_i\rho V_i^*)}
\end{equation}

\begin{lem} $h_V(W)\geq 0$, for every family $W_i$ of linear operators satisfying $\sum_i W_i^* W_i=I$.
\end{lem}
{\bf Proof} Note that, by definition,
$$h_V(W)=(\mathcal{U}h)(\rho_W)=\int_{\mathcal{M}_N} h d\mathcal{V}\delta_{\rho_W}$$
and the function $h$ (Shannon-Boltzmann entropy) is $\geq 0$. This proves the lemma. Another elementary proof is the following. Since $\rho_W$ is positive, we have that $\langle W_i\rho_W W_i^* v,v\rangle=\langle \rho_W W_i^* v,W_i^*v\rangle \geq 0$, $v\in\mathcal{H}_N$. So for $\{v_l\}_{l=1,\dots N}$ an orthonormal base for $\mathcal{H}_N$,
$$tr(W_i \rho_W W_i^*)=\sum_{l=1}^N\langle W_i \rho_W W_i^*v_l,v_l\rangle >0$$
Analogously the expression above holds for the $V_i\rho_W V_i^*$, and therefore also for $W_j V_i\rho_W V_i^* W_j^*$, because
$$\langle W_j V_i\rho_W V_i^* W_j^*v,v\rangle=\langle V_i\rho_W V_i^* W_j^*v,W_j^*v\rangle\geq 0$$
To conclude that $h_V(W)\geq 0$, we have to show that $tr(W_j V_i\rho_W V_i^*W_j^*)\leq tr(V_i\rho_W V_i^*)$. From $\sum_{i=1}^k W_i^*W_i=I$, we get
$$tr(W_j V_i\rho_W V_i^*W_j^*)=tr(W_j^*W_j V_i\rho_W V_i^*)\leq\sum_{j=1}^k tr(W_j^*W_j V_i\rho_W V_i^*)$$
$$=tr(\sum_{j=1}^k W_j^*W_j V_i\rho_W V_i^*)=tr(V_i\rho_W V_i^*)$$

\qed

{\bf Remark} For any fixed dynamics $V$, if we have that $W_m^*W_m=I$ for some $m$ then
the remaining $p_i$ must be zero, because of the condition $\sum_i
W_i^* W_i=I$. In this case we have $h_V(W)=0$. We also have that $h_V(W)\leq \log k$ and for any given dynamics $V$, $h_V(W)$ attains the maximum if we choose $W_i=1/\sqrt{k}I$, for each $i$, where $I$ denotes the identity operator.

\qee

Note that by the calculations made in section \ref{sec_calc_entr}, we have $h_V(W)=\mathcal{U}h(\rho_W)$, where $\mathcal{U}h(\rho)=\sum_i p_i(\rho)h(F_i(\rho))$.

\begin{lem}
Let $\mathcal{F}=(\mathcal{M}_N,F_i,p_i)$ be a QIFS, with $F_i$, $p_i$ in the form (\ref{din01}) and (\ref{prob01}). Suppose there is $\rho_0\in\mathcal{M}_N$ such that $\delta_{\rho_0}$ is the unique $\mathcal{V}$-invariant measure. Then $\Lambda_{\mathcal{F}}(\rho_0)=\rho_0$ ($\Lambda_{\mathcal{F}}$ is the operator associated to $\mathcal{F}$) and
$$\int \mathcal{U}^n h d\delta_{\rho_0}=\mathcal{U}^nh(\rho_0)=h(\rho_0),$$
for all $n\in\mathbb{N}$. Besides, $\mathcal{U}^nh(\rho_0)=\mathcal{U}h(\rho_0)$ and so
$$h_V(W)=\mathcal{U}^n h(\rho_0),$$
for all $n\in\mathbb{N}$.
\end{lem}
{\bf Proof} The fact that $\Lambda(\rho_0)=\rho_0$ follows from lemma \ref{lemabas1}, item 2. Also,
$$\mathcal{U}^nh(\rho_0)=\int \mathcal{U}^nh d\delta_{\rho_0}=\int h d\mathcal{V}^n\delta_{\rho_0}=\int h d\delta_{\rho_0}=h(\rho_0)$$
and
$$\mathcal{U}^nh(\rho_0)=\int \mathcal{U}^nh d\delta_{\rho_0}=\int h d\mathcal{V}^n\delta_{\rho_0}=\int h d\mathcal{V}\delta_{\rho_0}=\int \mathcal{U}h d\delta_{\rho_0}=\mathcal{U}h(\rho_0)$$

\qed

\begin{lem}
Let $\mu$ be a $\mathcal{V}$-invariant attractive measure. Then if $\rho_{\mu}$ is the barycenter of $\mu$ we have, for any $\rho$,
\begin{equation}\label{wcheck}
\lim_{n\to\infty}\mathcal{U}^nh(\rho)=\int \mathcal{U}h d\mu=\int h d\mu\leq h(\rho_{\mu})
\end{equation}
\end{lem}
{\bf Proof} The inequality follows from \cite{wsbook}, proposition 1.15. Also, by proposition \ref{pdualbas} we have
$$\lim_{n\to\infty}\mathcal{U}^nh(\rho)=\lim_{n\to\infty}\int\mathcal{U}^nh d\delta_{\rho}=\lim_{n\to\infty}\int \mathcal{U}h d\mathcal{V}^{n-1}\delta_{\rho}=\int \mathcal{U}h d\mu,$$
the last equality being true because of the weak convergence of $(\mathcal{V}^n\delta_{\rho})_{n\in \mathbb{N}}$. This proves the first equality in (\ref{wcheck}). Since $\int \mathcal{U}h d\mu=\int h d\mathcal{V}\mu=\int h d\mu$, we obtain the second equality.

\qed

\begin{lem}
Let $\mathcal{F}=(\mathcal{M}_N,F_i,p_i)$ be a QIFS, with $F_i$, $p_i$ in the form (\ref{din01}) and (\ref{prob01}).
Suppose that $\rho$ is the unique point such that $\Lambda_{\mathcal{F}}(\rho)=\rho$. Suppose that $F_i(\rho)=\rho$, $i=1,\dots, k$. Then
$$\mathcal{U}^nh(\rho)=h(\rho),$$
$n=1,2,\dots$, and therefore $h_V(W)$ does not depend on $n$.
\end{lem}
{\bf Proof} The proof follows by induction. Let $n=1$. We have:
$$\mathcal{U}h(\rho)=\sum_i p_i(\rho)h(F_i(\rho))=h(\rho)\sum_i p_i(\rho)=h(\rho)$$
And note that $\mathcal{U}^nh(\rho)=\mathcal{U}(\mathcal{U}^{n-1}h)(\rho)$, which concludes the proof.

\qed

\section{Entropy and Markov chains}\label{camarkov}

Let $V_i$, $W_i$ be linear operators, $i=1,\dots, k$, $\sum_{i=1}^k W_i^*W_i=I$. Suppose the $V_i$ are fixed and that they determine a dynamics given by $F_i:\mathcal{M}_N\to \mathcal{M}_N$, $i=1,\dots, k$. Define
\begin{equation}
P:=\{(p_1,\dots, p_k): p_i:\mathcal{M}_N\to\mathbb{R}^+,i=1,\dots, k ,\sum_{i=1}^k p_i(\rho)=1,\forall\rho\in\mathcal{M}_N\}
\end{equation}
$$
P':=P\cap\{(p_1,\dots, p_k):\exists W_i, i=1,\dots, k: p_i(\rho)=tr(W_i\rho W_i^*),$$
\begin{equation}
W_i \textrm{ linear }, \sum_i W_i^*W_i=I\}
\end{equation}

\begin{equation}
\mathcal{M}_F:=\{\mu\in M^1(\mathcal{M}_N): \exists p\in P' \textrm{ such that } \mathcal{V}_p\mu=\mu\},
\end{equation}
where $\mathcal{V}_p:M^1(\mathcal{M}_N)\to M^1(\mathcal{M}_N)$,
\begin{equation}
\mathcal{V}_p(\mu)(B):=\sum_{i=1}^k\int_{F_i^{-1}(B)}p_id\mu
\end{equation}

Note that a family $W:=\{W_i\}_{i=1,\dots, k}$ determines a QIFS $\mathcal{F}_{W}$,
$$\mathcal{F}_W=\{\mathcal{M}_N,F_i,p_i\}_{i=1,\dots, k}$$

Let $P=(p_{ij})_{i,j=1,\dots,N}$ be a stochastic, irreducible matrix. Let $p$ be the stationary vector of $P$. The entropy of $P$ is defined as
\begin{equation}\label{eestoc}
H(P):=-\sum_{i,j=1}^Np_i p_{ij}\log{p_{ij}}
\end{equation}
We consider a few examples which will be useful later in this work.

\begin{examp} (Homogeneous case, 4 matrices).
Let $N=2$, $k=4$ and
$$V_1=\left(
\begin{array}{cc}
\sqrt{p_{00}} & 0\\
0 & 0
\end{array}
\right),
\sp V_2=\left(
\begin{array}{cc}
0 & \sqrt{p_{01}}\\
0 & 0
\end{array}
\right)
,$$
$$V_3=\left(
\begin{array}{cc}
0 & 0\\
\sqrt{p_{10}} & 0
\end{array}
\right),
\sp V_4=\left(
\begin{array}{cc}
0 & 0\\
0 & \sqrt{p_{11}}
\end{array}
\right)
$$
Note that
$$\sum_i V_i^* V_i=\left(
\begin{array}{cc}
p_{00}+p_{10} & 0\\
0 & p_{01}+p_{11}
\end{array}
\right)$$
and so $\sum_i V_i^* V_i=I$ if we suppose that
$$P:=\left(
\begin{array}{cc}
p_{00} & p_{01}\\
p_{10} & p_{11}
\end{array}
\right)$$
is column-stochastic. We have
$$V_1\rho V_1^*=\left(
\begin{array}{cc}
p_{00}\rho_1 & 0\\
0 & 0
\end{array}
\right),\sp
V_2\rho V_2^*=\left(
\begin{array}{cc}
p_{01}\rho_4 & 0\\
0 & 0
\end{array}
\right)$$
$$V_3\rho V_3^*=\left(
\begin{array}{cc}
0 & 0\\
0 & p_{10}\rho_1
\end{array}
\right),\sp
V_4\rho V_4^*=\left(
\begin{array}{cc}
0 & 0\\
0 & p_{11}\rho_4
\end{array}
\right)
$$
so
$$tr(V_1\rho V_1^*)=p_{00}\rho_1,\sp tr(V_2\rho V_2^*)=p_{01}\rho_4$$
$$\sp tr(V_3\rho V_3^*)=p_{10}\rho_1,\sp tr(V_4\rho V_4^*)=p_{11}\rho_4$$

\bigskip
The fixed point of $\Lambda(\rho)=\sum_i V_i\rho V_i^*$ is
$$\rho_V=\left(
\begin{array}{cc}
\frac{p_{01}}{1-p_{00}+p_{01}} & 0\\
0 & \frac{1-p_{00}}{1-p_{00}+p_{01}}
\end{array}
\right)
$$

Let $\pi=(\pi_1,\pi_2)$ such that $P\pi=\pi$. We know that
\begin{equation}
\pi=(\frac{p_{01}}{1-p_{00}+p_{01}},\frac{1-p_{00}}{1-p_{00}+p_{01}})
\end{equation}
Then the nonzero entries of $\rho_V$ are the entries of $\pi$ and so we associate the fixed point of $P$ to the fixed point of a certain $\Lambda$ in a natural way. Let us calculate  $h_V(W)$. Note that $\Lambda$ defined above is associated to a homogeneous IFS. Then $W_i=V_i$, $i=1,\dots, k$ and
$$h_V(W)=h_V(V)$$
$$=-\sum_{i=1}^k \frac{tr(W_i\rho_V W_i^*)}{tr(V_i\rho_V V_i^*)}\sum_{j=1}^k tr\Big(W_j V_i\rho_V V_i^* W_j^*\Big)\log{\Big(\frac{tr(W_j V_i\rho_V V_i^* W_j^*)}{tr(V_i\rho_V V_i^*)}\Big)}$$
\begin{equation}
=-\sum_{i,j}tr\Big(V_j V_i\rho_V V_i^* V_j^*\Big)\log{\Big(\frac{tr(V_j V_i\rho_V  V_i^* V_j^*)}{tr(V_i\rho_V V_i^*)}\Big)}
\end{equation}
A simple calculation yields $H(P)=h_V(V)$, where $H(P)$ is the entropy of $P$, given by (\ref{eestoc}). This shows that the entropy of Markov chains is a particular case of the QIFS entropy.
\end{examp}

\qee

\begin{examp}\label{caso_util_a1}
(Nonhomogeneous case, 4 matrices). Let $N=2$, $k=4$ and
$$V_1=\left(
\begin{array}{cc}
\sqrt{p_{00}} & 0\\
0 & 0
\end{array}
\right),
\sp V_2=\left(
\begin{array}{cc}
0 & \sqrt{p_{01}}\\
0 & 0
\end{array}
\right)
$$
$$V_3=\left(
\begin{array}{cc}
0 & 0\\
\sqrt{p_{10}} & 0
\end{array}
\right),
\sp V_4=\left(
\begin{array}{cc}
0 & 0\\
0 & \sqrt{p_{11}}
\end{array}
\right)
$$
$$
W_1=\left(
\begin{array}{cc}
\sqrt{q_{00}} & 0 \\
0 & 0
\end{array}
\right),\sp W_2=\left(
\begin{array}{cc}
0 & \sqrt{q_{01}} \\
0 & 0
\end{array}
\right)$$
$$
W_3=\left(
\begin{array}{cc}
0 & 0\\
\sqrt{q_{10}} & 0
\end{array}
\right),\sp
W_4=\left(
\begin{array}{cc}
0 & 0\\
0 & \sqrt{q_{11}}
\end{array}
\right)$$

Note that
$$\sum_i V_i^* V_i=\left(
\begin{array}{cc}
p_{00}+p_{10} & 0\\
0 & p_{01}+p_{11}
\end{array}
\right),\sp \sum_i W_i^* W_i=\left(
\begin{array}{cc}
q_{00}+q_{10} & 0\\
0 & q_{01}+q_{11}
\end{array}
\right)$$
and so $\sum_i V_i^* V_i=\sum_i W_i^* W_i=I$ if we suppose that
$$P:=\left(
\begin{array}{cc}
p_{00} & p_{01}\\
p_{10} & p_{11}
\end{array}
\right),\sp Q:=\left(
\begin{array}{cc}
q_{00} & q_{01}\\
q_{10} & q_{11}
\end{array}
\right)\sp$$
are column-stochastic. Then
$$tr(V_1\rho V_1^*)=p_{00}\rho_1,\sp tr(V_2\rho V_2^*)=p_{01}\rho_4$$
$$\sp tr(V_3\rho V_3^*)=p_{10}\rho_1,\sp tr(V_4\rho V_4^*)=p_{11}\rho_4$$
$$tr(W_1\rho W_1^*)=q_{00}\rho_1,\sp tr(W_2\rho W_2^*)=q_{01}\rho_4$$
$$\sp tr(W_3\rho W_3^*)=q_{10}\rho_1,\sp tr(W_4\rho W_4^*)=q_{11}\rho_4$$

We want the fixed point of $\Lambda(\rho)=\sum_i p_i(\rho)F_i(\rho)$. This leads us to
$$\frac{q_{00}}{p_{00}}
\left(\begin{array}{cc}
p_{00}\rho_1 & 0\\
0 & 0
\end{array}
\right)+\frac{q_{01}}{p_{01}}
\left(\begin{array}{cc}
p_{01}\rho_4 & 0\\
0 & 0
\end{array}
\right)
+\frac{q_{10}}{p_{10}}
\left(\begin{array}{cc}
0 & 0\\
0 & p_{10}\rho_1
\end{array}
\right)
+\frac{q_{11}}{p_{11}}
\left(\begin{array}{cc}
0 & 0\\
0 & p_{11}\rho_4
\end{array}
\right)=\rho
$$
Note that the $p_{ij}$ cancel and so we obtain a calculation which is the same as the one obtained in the previous example. Hence
$$\rho_W=\left(
\begin{array}{cc}
\frac{q_{01}}{1-q_{00}+q_{01}} & 0\\
0 & \frac{1-q_{00}}{1-q_{00}+q_{01}}
\end{array}
\right),
$$
and its nonzero entries are the entries of the fixed point for the stochastic matrix $Q$. Calculating $h_V(W)$ gives
$$h_V(W)=-\sum_{i=1}^k \frac{tr(W_i\rho_W W_i^*)}{tr(V_i\rho_W V_i^*)}\sum_{j=1}^k tr\Big(W_j V_i\rho_W V_i^* W_j^*\Big)\log{\Big(\frac{tr(W_j V_i\rho_W  V_i^* W_j^*)}{tr(V_i\rho_W V_i^*)}\Big)}$$
\begin{equation}\label{entr_hvw_estoc}
=-\frac{q_{01}}{q_{01}+q_{10}}(q_{00}\log{q_{00}}+q_{10}\log{q_{10}})-\frac{q_{10}}{q_{01}+q_{10}}(q_{01}\log{q_{01}}+q_
{11}\log{q_{11}})=H(Q)
\end{equation}
So we have obtained a calculation which is analogous to the one for the homogeneous case. This result generalizes what we have seen in the previous example.
\end{examp}

\qee

\begin{examp}
(Homogeneous case, 2 matrices). Let $N=2$, $k=2$ and
$$V_1=\left(
\begin{array}{cc}
\sqrt{p_{00}} & 0\\
\sqrt{p_{10}} & 0
\end{array}
\right),
\sp V_2=\left(
\begin{array}{cc}
0 & \sqrt{p_{01}}\\
0 & \sqrt{p_{11}}
\end{array}
\right)
,$$
Note that, just as in the previous examples
$$\sum_i V_i^* V_i=\left(
\begin{array}{cc}
p_{00}+p_{10} & 0\\
0 & p_{01}+p_{11}
\end{array}
\right)$$
and so $\sum_i V_i^* V_i=I$ if we suppose
$$P:=\left(
\begin{array}{cc}
p_{00} & p_{01}\\
p_{10} & p_{11}
\end{array}
\right)$$
is column-stochastic. The fixed point for $\Lambda$ is
$$\rho_V=\left(
\begin{array}{cc}
\frac{p_{01}}{p_{01}+p_{10}} & \frac{p_{00}p_{10} p_{01}}{p_{01}+p_{10}}+\frac{p_{01}p_{11}p_{10}}{p_{01}+p_{10}}
\\
\frac{p_{00}p_{10} p_{01}}{p_{01}+p_{10}}+\frac{p_{01}p_{11}p_{10}}{p_{01}+p_{10}}
 & \frac{p_{10}}{p_{01}+p_{10}}
\end{array}
\right)
$$
The entries of the main diagonal of $\rho_V$ correspond to the entries of the fixed point of $P$. The entries of the secondary diagonal are a linear combination of the ones in the main diagonal. Then for the $V_i$ chosen we have
\begin{equation}
h_V(W)=h_V(V)=-\sum_{i,j}tr\Big(V_j V_i\rho_V V_i^* V_j^*\Big)\log{\Big(\frac{tr(V_j V_i\rho_V V_i^* V_j^*)}{tr(V_i\rho_V  V_i^*)}\Big)}=H(P)
\end{equation}
by an identical calculation made for the equation (\ref{entr_hvw_estoc}) from the previous example. In other words, the fact that the fixed point of $\Lambda$ is not diagonal does not change the calculations for the entropy.

\qee

\end{examp}

\begin{examp}
(Nonhomogeneous case, 2 matrices). Let $N=2$, $k=2$,
$$V_1=\left(
\begin{array}{cc}
\sqrt{p_{00}} & 0\\
\sqrt{p_{10}} & 0
\end{array}
\right),
\sp V_2=\left(
\begin{array}{cc}
0 & \sqrt{p_{01}}\\
0 & \sqrt{p_{11}}
\end{array}
\right)
$$
$$W_1=\left(
\begin{array}{cc}
\sqrt{q_{00}} & 0\\
\sqrt{q_{10}} & 0
\end{array}
\right),
\sp W_2=\left(
\begin{array}{cc}
0 & \sqrt{q_{01}}
\\
0 & \sqrt{q_{11}}
\end{array}
\right)
$$

As in the other examples, $\sum_i V_i^* V_i=\sum_i W_i^* W_i=I$ if we suppose
$$P:=\left(
\begin{array}{cc}
p_{00} & p_{01}\\
p_{10} & p_{11}
\end{array}
\right),\sp Q:=\left(
\begin{array}{cc}
q_{00} & q_{01}\\
q_{10} & q_{11}
\end{array}
\right)\sp$$
is column-stochastic. From
$$tr(V_1\rho V_1^*)=\rho_1,\sp tr(V_2\rho V_2^*)=\rho_4$$
$$tr(W_1\rho W_1^*)=\rho_1,\sp tr(W_2\rho W_2^*)=\rho_4$$
$$tr(W_1V_1\rho V_1^*W_1^*)=p_{00}\rho_1,\sp tr(W_2V_1\rho V_1^*W_2^*)=p_{10}\rho_1$$
$$tr(W_1V_2\rho V_2^*W_1^*)=p_{01}\rho_4,\sp tr(W_2V_2\rho V_2^*W_2^*)=p_{11}\rho_4$$
and a simple calculation, we get $h_V(W)=H(P)$.
\end{examp}

\qee

\begin{lem}
Let $V_{ij}$ be matrices of order $n$,
$$V_{ij}=\sqrt{p_{ij}}\vert i\rangle\langle j\vert$$
for $i,j=1,\dots, n$. Let $$\Lambda_P(\rho):=\sum_{i,j}V_{ij}\rho V_{ij}^*$$
where $P=(p_{ij})_{i,j=1,\dots, n}$. Then for all $n$, $\Lambda_P^n(\rho)=\Lambda_{P^n}(\rho)$.
\end{lem}
{\bf Proof} Note that
\begin{equation}
V_{kl}V_{ij}=\sqrt{p_{kl}}\sqrt{p_{ij}}\delta_{li}\vert k\rangle\langle j\vert
\end{equation}
so
$$\Lambda_P^2(\rho)=\Lambda_P(\sum_{i,j}V_{ij}\rho V_{ij}^*)=\sum_{k,l,i,j}V_{kl}V_{ij}\rho(V_{kl}V_{ij})^*$$
$$=\sum_{k,j}\sum_ip_{ki}p_{ij}\vert k\rangle\langle j\vert\rho\vert j\rangle\langle k\vert=\sum_{k,j}p_{kj}^2\vert k\rangle\langle j \vert\rho\vert j\rangle\langle k\vert=\Lambda_{P^2}(\rho)$$
The general case follows by iterating the above calculation.
\qed

\begin{cor}
Under the lemma hypothesis, we have $\lim_{n\to\infty}\Lambda_P^n(\rho)=\Lambda_{\pi}(\rho)$, where $\pi=\lim_{n\to\infty}P^n$ is the stochastic matrix which has all columns equal to the stationary vector for $P$.
\end{cor}

\section{Capacity-cost function and pressure}\label{sec_cc}

Recall that every trace preserving, completely positive (CP) mapping can be written in the Stinespring-Kraus form,
$$\Lambda(\rho)=\sum_{i=1}^k V_i\rho V_i^*,\sp \sum_{i=1}^k V_i^*V_i=I,$$
for $V_i$ linear operators. These mappings are also called {\bf quantum channels}. This is one of the main motivations for considering the class of operators (a generalization of the above ones) described in the present paper. These are natural objets in the analysis of certain problems in quantum computing.

\begin{defi} The {\bf Holevo capacity} for sending classic information via a quantum channel $\Lambda$ is defined as
\begin{equation}\label{hcap_01}
C_{\Lambda}:=\max_{\stackrel{p_i\in [0,1]}{\rho_i\in\mathcal{M}_N}}S\Big(\sum_{i=1}^n p_i\Lambda(\rho_i)\Big)-\sum_{i=1}^n p_i S\Big(\Lambda(\rho_i)\Big)
\end{equation}
where $S(\rho)=-tr(\rho\log\rho)$ is the von Neumann entropy. The maximum is, therefore, over all choices of $p_i$, $i=1,\dots, n$  and density operators $\rho_i$, for some $n\in \mathbb{N}$. The Holevo capacity establishes an upper bound on the amount of information that a quantum system contains \cite{nich}.
\end{defi}

\begin{defi} Let $\Lambda$ be a quantum channel. Define the {\bf minimum output entropy} as
$$H^{min}(\Lambda):=\min_{\vert\psi\rangle}S(\Lambda(\vert\psi\rangle\langle\psi\vert))$$
\end{defi}

{\bf Additivity conjecture} We have that
$$C_{\Lambda_1\otimes\Lambda_2}=C_{\Lambda_1}+C_{\Lambda_2}$$

{\bf Minimum output entropy conjecture} For any channels $\Lambda_1$ and $\Lambda_2$,
$$H^{min}(\Lambda_1\otimes\Lambda_2)=H^{min}(\Lambda_1)+H^{min}(\Lambda_2)$$

\bigskip
In \cite{shor}, is it shown that the additivity conjecture is equivalent to the minimum output entropy conjecture, and in \cite{hastings} a counterexample is obtained for this last conjecture.

\bigskip

 {\bf Remark}  Concerning QIFS, our interest in capacity is motivated by the following observation. Considering expression (\ref{hcap_01}), note that the term
\begin{equation}\label{qifs_e_cap}
\sum_{i=1}^n p_i S(\Lambda(\rho_i))
\end{equation}
is a convex combination of von Neumann entropies, in the same way as the QIFS entropy. So we see that given a QIFS, we can consider capacity functions, and the QIFS entropy arises in a natural way. For an example, we perform the following calculation. If $\lambda_i$ are the eigenvalues of $\rho$ then we can write
\begin{equation}\label{von_neumann_eigenvalues}
S(\rho)=-\sum_i \lambda_i\log\lambda_i
\end{equation}
Then write the QIFS entropy as
\begin{equation}\label{entropy_QIFS2}
h_V(W)=-\sum_{i,j=1}^k tr(W_i\rho_W W_i^*)
a_{ij}(\rho_W)\log a_{ij}(\rho_W)
\end{equation}
where
\begin{equation}\label{expressao_aij}
a_{ij}(\rho):=\frac{tr(W_jV_i\rho V_i^*W_j^* )}{tr(V_i\rho V_i^*)}
\end{equation}
We see that for $\rho_W\in\mathcal{M}_N$ and $i$ fixed, we have $\sum_i a_{ij}(\rho_W)=1$. Define for each $i$ the density operator
\begin{equation}\label{holevo1_later}
\rho_i:=\sum_j a_{ij}(\rho_W)\vert j\rangle\langle j\vert
\end{equation}
Then by (\ref{von_neumann_eigenvalues}),
\begin{equation}
S(\rho_i)=-\sum_j a_{ij}(\rho_W)\log a_{ij}(\rho_W)
\end{equation}
By (\ref{entropy_QIFS2}), we can write
\begin{equation}\label{holevo2_later}
h_V(W)(\rho_W)=\sum_{i=1}^k tr(W_i\rho_W W_i^*)S(\rho_i)
\end{equation}
A {\bf Positive Operator-Valued Measurement} (POVM) is described by a set of positive operators $P_i$ (POVM elements) such that $\sum_i P_i=I$. If the measurement is performed on a system described by the state vector $\vert\psi\rangle$, then the probability of obtaining $i$ as the outcome is given by
\begin{equation}
p_i=\langle\psi\vert F_i\vert\psi\rangle
\end{equation}
Note that a QIFS $\mathcal{F}$ induced by linear $V_i$ and $W_i$, contains a POVM by taking $W_i^*W_i$ as POVM elements. If $X$ is a random variable that takes values $p_1,\dots p_k$ then the Shannon entropy is $H(X)=-\sum_i p_i\log p_i$ and the joint entropy of variables $X$ and $Y$ is
\begin{equation}
H(X,Y):=-\sum_{x,y} p(x,y)\log p(x,y)
\end{equation}
where $p(x,y)$ is the probability that $X=x$ and $X=y$. The mutual information $I(X:Y)$ is defined by $I(X:Y):=H(X)+H(Y)-H(X,Y)$. Then, considering the QIFS entropy we can state the Holevo bound in the following way: first consider a QIFS $\mathcal{F}$ such that there is a unique attractive measure which is invariant for the Markov operator $\mathcal{V}$ associated to $\mathcal{F}$. Let $\rho_W$ be the barycenter of such measure.

\begin{teo} {\bf (Holevo bound for QIFS)} Suppose $\mathcal{F}$ is induced by linear operators $V_i$ and $W_i$ with $\sum_i W_i^*W_i=I$ and for each $i=1,\dots, k$ write $p_i=tr(W_i\rho_W W_i^*)$ and $\rho_i=\sum_j a_{ij}(\rho_W)\vert j\rangle\langle j\vert$, where $a_{ij}$ is given by (\ref{expressao_aij}). Suppose Alice prepares a mixed state $\rho_X$ chosen from the ensemble $\{\rho_1,\dots, \rho_k\}$ with probabilities $\{p_1,\dots, p_k\}$ (that is, we assume $\rho_X$ is a state determined by a random variable $X$ such that it assumes the value $\rho_i$ with probability $p_i$). Suppose Bob performs a POVM measurement on that state with POVM elements $\{P_i\}_{i=1,\dots, m}$ and measurement outcome described by a random variable $Y$. Then, by writing $\rho=\sum_i p_i\rho_i$, we have
\begin{equation}\label{holevo_qifs}
I(X:Y)\leq S(\rho)-h_V(W)(\rho_W)=:\xi(\mathcal{E})
\end{equation}
\end{teo}
The number $\xi(\mathcal{E})$ is the Holevo information of the ensemble given by $\mathcal{E}=\{\rho_i; p_i\}_{i=1,\dots, k}$. We see that (\ref{holevo_qifs}) holds by applying the Holevo bound for the von Neumann entropy (see \cite{nich}) together with (\ref{holevo1_later}) and (\ref{holevo2_later}).

\qee

We are also interested in a different class of problems which concern maximization (and not minimization) of entropy plus a given potential (a cost) \cite{gray},\cite{hayashi1},\cite{hayashi2}.

\begin{defi}
Let $M_F$ be the set of invariant measures defined in the section \ref{camarkov} and let $H$ be a hermitian operator. For $\mu\in\mathcal{M}_F$ let $\rho_{\mu}$ be its barycenter. Define the capacity-cost function $C:\mathbb{R}^+\to\mathbb{R}^+$ as
\begin{equation}
C(a):=\max_{\mu\in\mathcal{M}_F} \{ h_{W,V}(\rho_{\mu}): tr(H\rho_{\mu})\leq a\}
\end{equation}
\end{defi}
The following analysis is inspired in \cite{lopescra}. There is a relation between the cost-capacity function and the variational problem for pressure. In fact, let $F:\mathbb{R}^+\to\mathbb{R}^+$ be the function given by
\begin{equation}
F(\lambda):=\sup_{\mu\in\mathcal{M}_F}\{h_{W,V}(\rho_{\mu})-\lambda tr(H\rho_{\mu})\}
\end{equation}
We have the following fact. There is a unique probability measure $\nu_0\in\mathcal{M}_F$ such that
$$F(\lambda)=h_{W,V}(\rho_{\nu_0})-\lambda tr(H\rho_{\nu_0})$$
Also, we have the following lemma:
\begin{lem}
Let $\lambda\leq 0$, and $\hat{a}=tr(H\rho_{\nu_0})$. Then
\begin{equation}
C(\hat{a})=h_{W,V}(\rho_{\nu_0})
\end{equation}
\end{lem}
{\bf Proof} Let $\nu\in\mathcal{M}_F$, $\nu\neq\nu_0$, with $tr(H\rho_\nu)\leq \hat{a}=tr(H\rho_{\nu_0})$. Then
$$h_{W,V}(\rho_{\nu})-\lambda tr(H\rho_\nu)<h_{W,V}(\rho_{\nu_0})-\lambda tr(H\rho_{\nu_0})$$
so
$$h_{W,V}(\rho_{\nu})<h_{W,V}(\rho_{\nu_0})$$
Hence
$$h_{W,V}(\rho_{\nu_0})=\sup_{\mu\in\mathcal{M}_F}\{h_{W,V}(\rho_\mu):tr(H\rho_\mu)\leq\hat{a}\}=C(\hat{a})$$

\qed

\section{Analysis of the pressure problem}\label{sec_analysis1}

Let $V_i$, $W_i$, $H_i$ be linear operators, $i=1,\dots, k$, with $\sum_i W_i^*W_i=I$ and let
\begin{equation}\label{um_pot_int}
H\rho:=\sum_{i=1}^k H_i\rho H_i^*
\end{equation}
a hermitian operator. We are interested in obtaining a version of the variational principle of pressure for our context. We will see that the pressure will be a maximum whenever we have a certain relation between the potential $H$ and the probability distribution considered (represented here by the $W_i$). We begin by fixing a dynamics, given by the $V_i$. From the reasoning described below, it will be natural to consider as definition of pressure the maximization among the possible stationary $W_i$ of the expression
$$
h_V(W)+ \sum_{j=1}^k \log \Big( tr(H_j\rho_{\beta}  H_j^*)tr(V_j\rho_{\beta} V_j^*)\Big) tr(W_j \rho_W  W_j^*)$$
where $\rho_\beta$ is the eigenstate of a certain Ruelle operator, described below. We begin our analysis by using the following elementary lemma.

\begin{lem}\label{lemalog2}
\cite{par} If $r_1,\dots , r_k$ and $q_1,\dots ,q_k$ are two probability distributions over $1,\dots, k$,   such that $r_j>0$, $j=1,\dots, k$, then
\begin{equation}
-\sum_{j=1}^kq_j\log{q_j}+\sum_{j=1}^k q_j\log{r_j}\leq 0
\end{equation}
and equality holds if and only if $r_j=q_j$, $j=1,\dots, k$.
\end{lem}

\bigskip

The potential given by (\ref{um_pot_int}), together with the $V_i$, induces an operator given by
\begin{equation}\label{operador_usa1}
\mathcal{L}_H(\rho):=\sum_{i=1}^k tr(H_i\rho H_i^*)V_i\rho V_i^*
\end{equation}
By proposition \ref{proavl1} we know that such operator admits an eigenvalue $\beta$ with its associated eigenstate $\rho_{\beta} $. Then $\mathcal{L}_H(\rho_{\beta} )=\beta\rho_{\beta} $ implies
\begin{equation}\label{soma_1_a}
\sum_{i=1}^k tr(H_i\rho_{\beta}  H_i^*)V_i\rho_{\beta}  V_i^*=\beta\rho_{\beta}
\end{equation}
In coordinates, (\ref{soma_1_a}) can be written as
\begin{equation}\label{soma_1_b}
\sum_{i=1}^k tr(H_i\rho_\beta H_i^*)(V_i\rho_\beta V_i^*)_{jl}=\beta(\rho_\beta)_{jl}
\end{equation}

\bigskip

{\bf Remark} Comparing the above calculation with the problem of finding an eigenvalue $\lambda$ of a matrix $A=(a_{ij})$, we have that equation (\ref{soma_1_a}) can be seen as the analogous of the expression
\begin{equation}
lE^A=\lambda l
\end{equation}
Above, the matrix $A$ plays the role of a potential, $E^A$ denotes the matrix with entries $e^{a_{ij}}$ and $l_j$ denotes the $j$-th coordinate of the left eigenvector $l$ associated to the eigenvalue $\lambda$. In coordinates,
\begin{equation}
\sum_i l_i e^{a_{ij}}=\lambda l_j, \sp i,j=1,\dots, k
\end{equation}
\qee

From this point we can perform two calculations. First,
considering (\ref{soma_1_a}) we will take the trace of such
equation in order to obtain a scalar equation. In spite of the
fact that taking the trace makes us lose part of the information
given by the eigenvector equation, we are still able to obtain a
version of what we will call the {\bf basic inequality}, which can
be seen as a QIFS version of the variational principle of
pressure. However, there is an algebraic drawback to this
approach, namely, that we will not be able to recover the classic
variational problem as a particular case of such inequality (such
disadvantage is a consequence of taking the trace, clearly). The
second calculation begins at equation (\ref{soma_1_b}), the
coordinate equations associated to the matrix equation for the
eigenvectors. In this case we also obtain a basic inequality, but
then we will have the classic variational problem of pressure as a
particular case.

\bigskip

An important question which is of our interest, regarding both calculations mentioned above, is to ask whether it is possible for a given system to attain its maximum pressure. It is not clear that given any dynamics, we can obtain a measure reaching such a maximum. With respect to our context, we will remark a natural condition on the dynamics which allows us to determine expressions for the measure which maximizes the pressure. Now we perform the calculations mentioned above.

\bigskip
Based on (\ref{soma_1_a}), define
\begin{equation}
r_j=\frac{1}{\beta}tr(H_j\rho_{\beta}  H_j^*)tr(V_j\rho_{\beta}  V_j^*)
\end{equation}
So we have $\sum_j r_j=1$. Let
\begin{equation}
q_j^i:=tr\Big(\frac{W_j V_i\rho_W V_i^* W_j^{*}}{tr(V_i\rho_W V_i^*)}\Big)
\end{equation}
where, as before, $\rho_W$ is the fixed point associated to the operator $\Lambda_{\mathcal{F}_W}$
\begin{equation}
\Lambda_{\mathcal{F}_W}(\rho):=\sum_{i=1}^k p_i(\rho)F_i(\rho)
\end{equation}
induced by the QIFS $(\mathcal{M}_N,F_i,p_i)_{i=1,\dots, k}$,
$$F_i(\rho)=\frac{V_i\rho V_i^*}{tr(V_i\rho V_i^*)}$$
and
$$p_i(\rho)=tr(W_i\rho W_i^*)$$
Note that we have
$$\sum_{j=1}^kq_j^i=\frac{1}{tr(V_i\rho_W V_i^*)}\sum_{j=1}^k tr(W_j^{*}W_j V_i\rho_W V_i^*)$$
$$=\frac{1}{tr(V_i\rho_W V_i^*)}tr(\sum_{j=1}^k W_j^{*}W_j V_i\rho_W V_i^*) =1$$

Then we can apply lemma \ref{lemalog2} for $r_j$, $q_j^i$, $j=1,\dots k$, with $i$ fixed, to obtain
$$-\sum_j tr\Big(\frac{W_j V_i\rho_W V_i^* W_j^{*}}{tr(V_i\rho_W V_i^*)}\Big)\log tr\Big(\frac{W_j V_i\rho_W V_i^* W_j^{*}}{tr(V_i\rho_W V_i^*)}\Big)$$
\begin{equation}
+\sum_j tr\Big(\frac{W_j V_i\rho_W V_i^* W_j^{*}}{tr(V_i\rho_W V_i^*)}\Big)\log \Big(\frac{1}{\beta}tr(H_j\rho_{\beta}  H_j^*)tr(V_j\rho_{\beta}  V_j^*)\Big)\leq 0
\end{equation}
and equality holds if and only if for all $i,j$,
\begin{equation}
\frac{1}{\beta}tr(H_j\rho_{\beta}  H_j^*)tr(V_j\rho_{\beta}  V_j^*)=\frac{tr(W_j V_i\rho_W V_i^* W_j^{*})}{tr(V_i\rho_W V_i^*)}
\end{equation}
Then
$$-\sum_j tr\Big(\frac{W_j V_i\rho_W V_i^* W_j^{*}}{tr(V_i\rho_W V_i^*)}\Big)\log tr\Big(\frac{W_j V_i\rho_W V_i^* W_j^{*}}{tr(V_i\rho_W V_i^*)}\Big)$$
$$+\sum_j tr\Big(\frac{W_j V_i\rho_W V_i^* W_j^{*}}{tr(V_i\rho_W V_i^*)}\Big)\log \Big( tr(H_j\rho_{\beta}  H_j^*)tr(V_j\rho_{\beta}  V_j^*)\Big)$$
$$\leq \sum_j tr\Big(\frac{W_j V_i\rho_W V_i^* W_j^{*}}{tr(V_i\rho_W V_i^*)}\Big)\log\beta
$$
which is equivalent to
$$-\sum_j tr\Big(\frac{W_j V_i\rho_W V_i^* W_j^{*}}{tr(V_i\rho_W V_i^*)}\Big)\log tr\Big(\frac{W_j V_i\rho_W V_i^* W_j^{*}}{tr(V_i\rho_W V_i^*)}\Big)$$
\begin{equation}
+\sum_j \frac{tr(W_j V_i\rho_W V_i^* W_j^{*})}{tr(V_i\rho_W V_i^*)}\log \Big( tr(H_j\rho_{\beta}  H_j^*)tr(V_j\rho_{\beta}  V_j^*)\Big)\leq \log\beta
\end{equation}
Multiplying by $tr(W_i\rho_W W_i^*)$ and summing over the $i$ index, we have
$$h_V(W)+\sum_j \log \Big( tr(H_j\rho_{\beta}  H_j^*)tr(V_j\rho_{\beta}  V_j^*)\Big) \sum_i \frac{tr(W_i\rho_W W_i^*)}{tr(V_i\rho_W V_i^*)} tr(W_j V_i\rho_W V_i^* W_j^{*})$$
\begin{equation}\label{pv_pr01xx}
\leq \sum_i tr(W_i\rho_W W_i^*)\log\beta=\log\beta
\end{equation}
and equality holds if and only if for all $i,j$,
\begin{equation}
\frac{1}{\beta}tr(H_j\rho_{\beta}  H_j^*)tr(V_j\rho_{\beta}  V_j^*)=\frac{tr(W_j V_i\rho_W V_i^* W_j^{*})}{tr(V_i\rho_W V_i^*)}
\end{equation}

Let us rewrite inequality (\ref{pv_pr01xx}). First we use the fact that $\rho_W$ is a fixed point of $\Lambda_{\mathcal{F}_W}$,
\begin{equation}
\sum_{i=1}^k tr(W_i\rho_W W_i^*)\frac{V_i\rho_W V_i^*}{tr(V_i\rho_W V_i^*)}=\rho_W
\end{equation}
Now we compose both sides of the equality above with the operator
\begin{equation}
\sum_{j=1}^k \log \Big( tr(H_j\rho_{\beta}  H_j^*)tr(V_j\rho_{\beta}  V_j^*)\Big) W_j^*W_j
\end{equation}
and then we obtain
$$\sum_{i=1}^k tr(W_i\rho_W W_i^*)\frac{V_i\rho_W V_i^*}{tr(V_i\rho_W V_i^*)}\sum_{j=1}^k \log \Big( tr(H_j\rho_{\beta}  H_j^*)tr(V_j\rho_{\beta}  V_j^*)\Big) W_j^*W_j$$
\begin{equation}
=\rho_W \sum_{j=1}^k \log \Big( tr(H_j\rho_{\beta}  H_j^*)tr(V_j\rho_{\beta}  V_j^*)\Big) W_j^*W_j
\end{equation}
Reordering terms we get
$$\sum_{j=1}^k \log \Big( tr(H_j\rho_{\beta}  H_j^*)tr(V_j\rho_{\beta}  V_j^*)\Big) \sum_{i=1}^k \frac{tr(W_i\rho_W W_i^*)}{tr(V_i\rho_W V_i^*)} V_i\rho_W V_i^* W_j^*W_j$$
\begin{equation}
=\rho_W \sum_{j=1}^k \log \Big( tr(H_j\rho_{\beta}  H_j^*)tr(V_j\rho_{\beta}  V_j^*)\Big) W_j^*W_j
\end{equation}
Taking the trace on both sides we get
$$\sum_{j=1}^k \log \Big( tr(H_j\rho_{\beta}  H_j^*)tr(V_j\rho_{\beta}  V_j^*)\Big) \sum_{i=1}^k \frac{tr(W_i\rho_W W_i^*)}{tr(V_i\rho_W V_i^*)} tr(W_j V_i\rho_W V_i^* W_j^*)$$
\begin{equation}\label{subst_aquixx}
=\sum_{j=1}^k \log \Big( tr(H_j\rho_{\beta}  H_j^*)tr(V_j\rho_{\beta} V_j^*)\Big) tr(\rho_W  W_j^*W_j)
\end{equation}
Note that the left hand side of (\ref{subst_aquixx}) is one of the sums appearing in (\ref{pv_pr01xx}). Therefore replacing (\ref{subst_aquixx}) into (\ref{pv_pr01xx}) gives our main result.
\begin{teo} Let $\mathcal{F}_W$ be a QIFS such that there is a unique attractive invariant measure for the associated Markov operator $\mathcal{V}$. Let $\rho_W$ be the barycenter of such measure and let $\rho_\beta$ be an eigenstate of $\mathcal{L}_H(\rho)$ with eigenvalue $\beta$. Then
\begin{equation}\label{p_var_geralxx} 
h_V(W)+ \sum_{j=1}^k \log \Big( tr(H_j\rho_{\beta}  H_j^*)tr(V_j\rho_{\beta} V_j^*)\Big) tr(W_j \rho_W  W_j^*)
 \leq \log\beta
\end{equation}
and equality holds if and only if for all $i,j$,
\begin{equation}\label{cond_igualdadexx}
\frac{1}{\beta}tr(H_j\rho_{\beta}  H_j^*)tr(V_j\rho_{\beta}  V_j^*)=\frac{tr(W_j V_i\rho_W V_i^* W_j^{*})}{tr(V_i\rho_W V_i^*)}
\end{equation}
\end{teo}
In section \ref{sec_analysis4} we make some considerations about certain cases in which we can reach an equality in (\ref{p_var_geralxx}).

\qee

For the calculations regarding expression (\ref{soma_1_b}), define
\begin{equation}
r_{jlm}=\frac{1}{\beta}tr(H_j\rho_{\beta}  H_j^*)\frac{(V_j\rho_{\beta}  V_j^*)_{lm}}{(\rho_\beta)_{lm}}
\end{equation}
Then we have $\sum_j r_{jlm}=1$. Let
\begin{equation}
q_{ij}:=tr\Big(\frac{W_j V_i\rho_W V_i^* W_j^{*}}{tr(V_i\rho_W V_i^*)}\Big)
\end{equation}
A calculation similar to the one we have made for (\ref{p_var_geralxx}) gives us
$$h_V(W)+ \sum_{j=1}^k tr(W_j \rho_W  W_j^*) \log tr(H_j\rho_{\beta}  H_j^*)$$
\begin{equation}
+\sum_{j=1}^k tr(W_j \rho_W  W_j^*) \log{\Big(\frac{(V_j\rho_{\beta} V_j^*)_{lm}}{(\rho_\beta)_{lm}}\Big)}
 \leq \log\beta
\end{equation}
and equality holds if and only if for all $i,j,l,m$,
\begin{equation}
\frac{1}{\beta}tr(H_j\rho_{\beta}  H_j^*)\frac{(V_j\rho_{\beta} V_j^*)_{lm}}{(\rho_\beta)_{lm}}=\frac{tr(W_j V_i\rho_W V_i^* W_j^{*})}{tr(V_i\rho_W V_i^*)}
\end{equation}

\qee

\section{Revisiting the eigenvalue problem}\label{sec_analysis2}

Consider the operator
\begin{equation}
\mathcal{L}_H(\rho)=\sum_{i=1}^k tr(H_i\rho H_i^*)V_i\rho V_i^*
\end{equation}
induced by a fixed dynamics $V_i$ $i=1,\dots, k$, $V_i$ linear, and by $H\rho:=\sum_i H_i\rho H_i^*$, $H_i$ linear. The eigenvalues equation for $\mathcal{L}_H$ written in coordinates gives us the following system, for $k=2$:

$$tr(H_1\rho_\beta H_1^*)(v_{11}^2\rho_{11}+2v_{11}v_{12}\rho_{12}+v_{12}^2\rho_{22})$$
\begin{equation}\label{cl_s1}
+tr(H_2\rho_\beta H_2^*)(w_{11}^2\rho_{11}+2w_{11}w_{12}\rho_{12}+w_{12}^2\rho_{22})=\beta\rho_{11}
\end{equation}

$$tr(H_1\rho_\beta H_1^*)(v_{21}v_{11}\rho_{11}+(v_{21}v_{12}+v_{22}v_{11})\rho_{12}+v_{22}v_{12}\rho_{22})$$
\begin{equation}\label{cl_s2}
+tr(H_2\rho_\beta H_2^*)(w_{21}w_{11}\rho_{11}+(w_{21}w_{12}+w_{22}w_{11})\rho_{12}+w_{22}w_{12}\rho_{22})=\beta\rho_{12}
\end{equation}

$$tr(H_1\rho_\beta H_1^*)(v_{21}^2\rho_{11}+2v_{21}v_{22}\rho_{12}+v_{22}^2\rho_{22})$$
\begin{equation}\label{cl_s3}
+tr(H_2\rho_\beta H_2^*)(w_{21}^2\rho_{11}+2w_{21}w_{22}\rho_{12}+w_{22}^2\rho_{22})=\beta\rho_{22}
\end{equation}

And we can also write, for $i=1,2$,
\begin{equation}
tr(H_i\rho H_i^*)=((h_{11}^i)^2+(h_{12}^i)^2)\rho_{11}+2(h_{11}^ih_{12}^i+h_{12}^ih_{22}^i)\rho_{12}+((h_{12}^i)^2+(h_{22}^i)^2)\rho_{22}
\end{equation}

\qee

Fix $H_1$, $H_2$, let $V_1$, $V_2$ be defined by
\begin{equation}
V_1=\left(
\begin{array}{cc}
v_{11} & v_{12} \\
0 & 0
\end{array}
\right),\sp V_2=\left(
\begin{array}{cc}
0 & 0 \\
w_{21} & w_{22}
\end{array}
\right)
\end{equation}
then we have, by (\ref{cl_s1})-(\ref{cl_s3}) that $\rho_{12}=0$ and
\begin{equation}
tr(H_1\rho_\beta H_1^*)(v_{11}^2\rho_{11}+v_{12}^2\rho_{22})=\beta\rho_{11}
\end{equation}
\begin{equation}
tr(H_2\rho_\beta H_2^*)(w_{21}^2\rho_{11}+w_{22}^2\rho_{22})=\beta\rho_{22}
\end{equation}
that is,
\begin{equation}
[((h_{11}^1)^2+(h_{12}^1)^2)\rho_{11}+((h_{12}^1)^2+(h_{22}^1)^2)\rho_{22}](v_{11}^2\rho_{11}+v_{12}^2\rho_{22})=\beta\rho_{11}
\end{equation}
\begin{equation}
[((h_{11}^2)^2+(h_{12}^2)^2)\rho_{11}+((h_{12}^2)^2+(h_{22}^2)^2)\rho_{22}](w_{21}^2\rho_{11}+w_{22}^2\rho_{22})=\beta\rho_{22}
\end{equation}

Also, suppose that
\begin{equation}
v_{11}=v_{12}=w_{21}=w_{22}=1
\end{equation}

Then we get
\begin{equation}\label{syst-1}
((h_{11}^1)^2+(h_{12}^1)^2)\rho_{11}+((h_{12}^1)^2+(h_{22}^1)^2)\rho_{22}=\beta\rho_{11}
\end{equation}
\begin{equation}\label{syst-2}
((h_{11}^2)^2+(h_{12}^2)^2)\rho_{11}+((h_{12}^2)^2+(h_{22}^2)^2)\rho_{22}=\beta\rho_{22}
\end{equation}

\bigskip

Let $A=(a_{ij})$ be a matrix with positive entries and consider the problem of finding its eigenvalues and eigenvectors. Then from
\begin{equation}\label{syst01}
a_{11}v_1+a_{12}v_2=\beta v_1
\end{equation}
\begin{equation}\label{syst02}
a_{21}v_1+a_{22}v_2=\beta v_2
\end{equation}
we see that the systems (\ref{syst-1})-(\ref{syst-2}) and (\ref{syst01})-(\ref{syst02}) are the same if we choose
\begin{equation}
a_{11}=(h_{11}^1)^2+(h_{12}^1)^2,\sp a_{12}=(h_{12}^1)^2+(h_{22}^1)^2
\end{equation}
\begin{equation}
a_{21}=(h_{11}^2)^2+(h_{12}^2)^2,\sp a_{22}=(h_{12}^2)^2+(h_{22}^2)^2
\end{equation}
We conclude that Perron's classic eigenvalue problem is a particular case of the problem associated to $\mathcal{L}_H$ acting on matrices. In fact, if we fix
\begin{equation}\label{exmm01}
V_1=\left(
\begin{array}{cc}
1 & 1 \\
0 & 0
\end{array}
\right),\sp V_2=\left(
\begin{array}{cc}
0 & 0 \\
1 & 1
\end{array}
\right)
\end{equation}
and given $A$ a matrix with positive entries, choose
\begin{equation}\label{exmm02}
H_1=\left(
\begin{array}{cc}
\sqrt{a_{11}} & 0 \\
0 & \sqrt{a_{12}}
\end{array}
\right),\sp H_2=\left(
\begin{array}{cc}
\sqrt{a_{21}} & 0 \\
0 & \sqrt{a_{22}}
\end{array}
\right)
\end{equation}
Then the operator $\mathcal{L}_H$ has a diagonal eigenstate
\begin{equation}
\rho_\beta=\left(
\begin{array}{cc}
\rho_{11} & 0 \\
0 & \rho_{22}
\end{array}
\right)
\end{equation}
associated to the eigenvalue $\beta$, and we have that, defining $v=(\rho_{11},\rho_{22})$, we get $Av=\beta v$.

\begin{examp}
Let
\begin{equation}
V_1=\left(
\begin{array}{cc}
1 & 1 \\
0 & 0
\end{array}
\right),\sp V_2=\left(
\begin{array}{cc}
0 & 0 \\
1 & 1
\end{array}
\right),\sp A=\left(
\begin{array}{cc}
1 & 4 \\
3 & \frac{1}{2}
\end{array}
\right)
\end{equation}
Then $Av=\beta v$ leads us to
\begin{equation}\label{mesmo_s1}
v_1+4v_2=\beta v_1
\end{equation}
\begin{equation}\label{mesmo_s2}
3v_1+\frac{1}{2}v_2=\beta v_2
\end{equation}
The eigenvalues are
$$\frac{3}{4}\pm\frac{1}{4}\sqrt{193}$$
with eigenvectors
$$\frac{1}{1\pm\frac{1}{12}+\frac{1}{12}\sqrt{193}}(\frac{1}{12}\pm\frac{1}{12}\sqrt{193},1)$$
Then we have $\beta=\frac{3}{4}+\frac{1}{4}\sqrt{193}$, $v=\frac{1}{1+\frac{1}{12}+\frac{1}{12}\sqrt{193}}(\frac{1}{12}+\frac{1}{12}\sqrt{193},1)$ such that $Av=\beta v$. Let
\begin{equation}
H_1=\left(
\begin{array}{cc}
\sqrt{a_{11}} & 0 \\
0 & \sqrt{a_{12}}
\end{array}
\right)=\left(
\begin{array}{cc}
1 & 0 \\
0 & 2
\end{array}
\right),\sp H_2=\left(
\begin{array}{cc}
\sqrt{a_{21}} & 0 \\
0 & \sqrt{a_{22}}
\end{array}
\right)=\left(
\begin{array}{cc}
\sqrt{3} & 0 \\
0 & \frac{1}{\sqrt{2}}
\end{array}
\right)
\end{equation}
Then solving $\mathcal{L}_H(\rho)=\beta\rho$ gives us $\rho_{12}=0$ and
\begin{equation}
\rho_{11}+4\rho_{22}=\beta\rho_{11}
\end{equation}
\begin{equation}
3\rho_{11}+\frac{1}{2}\rho_{22}=\beta\rho_{22}
\end{equation}
which is the same system as (\ref{mesmo_s1})-(\ref{mesmo_s2}). So
$\beta=\frac{3}{4}+\frac{1}{4}\sqrt{193}$ and the corresponding eigenstate, since $\rho_{12}=0$, is
\begin{equation}
\rho=\left(
\begin{array}{cc}
\frac{\frac{1}{12}+\frac{1}{12}\sqrt{193}}{1+\frac{1}{12}+\frac{1}{12}\sqrt{193}} & 0 \\
0 & \frac{1}{1+\frac{1}{12}+\frac{1}{12}\sqrt{193}}
\end{array}
\right)
\end{equation}
\end{examp}

\qee

\section{Some classic inequality calculations}\label{sec_analysis3}

A natural question is to ask whether the maximum among normalized
$W_i$, $i=1,\dots,k,$ for the pressure problem associated to a
given potential is realized as the logarithm of the main
eigenvalue of a certain Ruelle operator associated to the
potential $H_i$, $i=1,\dots ,k.$ This problem will be considered
in this section and also in the next one.

We begin by recalling a classic inequality. Consider
\begin{equation}\label{adbas1}
-\sum_{j=1}^kq_j\log{q_j}+\sum_{j=1}^k q_j\log{r_j}\leq 0
\end{equation}
given by lemma \ref{lemalog2}. Let $A$ be a matrix. If $v$ denotes the left eigenvector of matrix $E^A$ (such that each entry is $e^{a_{ij}}$), then $vE^A=\beta v$ can be written as
\begin{equation}
\sum_i v_{i}e^{a_{ij}}=\beta v_{j},\sp \forall j
\end{equation}
Define
\begin{equation}
r_{ij}:=\frac{e^{a_{ij}}v_{i}}{\beta v_j}
\end{equation}
So $\sum_i r_{ij}=1$. Let $q_{ij}>0$ such that $\sum_i q_{ij}=1$. By (\ref{adbas1}), we have
\begin{equation}
-\sum_{i=1}^kq_{ij}\log{q_{ij}}+\sum_{i=1}^k q_{ij}\log{\frac{e^{a_{ij}}v_{i}}{\beta v_j}}\leq 0
\end{equation}
That is,
\begin{equation}
-\sum_{i=1}^kq_{ij}\log{q_{ij}}+\sum_{i=1}^k q_{ij}a_{ij}+\sum_{i=1}^k q_{ij}(\log{v_{i}}-\log{v_{j}})\leq \log{\beta}
\end{equation}
Let $Q$ be a matrix with entries $q_{ij}$, let $\pi=(\pi_1,\dots, \pi_k)$ be the stationary vector associated to $Q$. Since $\sum_i q_{ij}=1$, $Q$ is column-stochastic so we write $Q\pi=\pi$. Multiplying the above inequality by $\pi_j$ and summing the $j$ index, we get
\begin{equation}
-\sum_j\pi_j\sum_i q_{ij}\log{q_{ij}}+\sum_j\pi_j\sum_i q_{ij}a_{ij}+\sum_j\pi_j\sum_i q_{ij}(\log{v_{i}}-\log{v_{j}})\leq \log{\beta}
\end{equation}
In coordinates, $Q\pi=\pi$ is $\sum_j q_{ij}\pi_j =\pi_i$, for all $i$. Then
$$-\sum_j\pi_j\sum_i q_{ij}\log{q_{ij}}+\sum_j\pi_j\sum_i q_{ij}a_{ij}$$
\begin{equation}
+\sum_j\pi_j\sum_i q_{ij}\log{v_{i}} -\sum_j\pi_j\sum_i q_{ij}\log{v_{j}}\leq \log{\beta}
\end{equation}
These calculations are well-known and gives us the following
inequality:
\begin{equation}\label{desig_classica}
-\sum_j\pi_j\sum_i q_{ij}\log{q_{ij}}+\sum_j\pi_j\sum_i q_{ij}a_{ij}\leq \log{\beta}
\end{equation}

\begin{defi} We call inequality (\ref{desig_classica}) the {\bf classic inequality} associated to the matrix $A$ with positive entries, and stochastic matrix $Q$.
\end{defi}

\begin{defi} For fixed $k$, and $l,m=1,\dots, k$ we call the inequality
$$h_V(W)+ \sum_{j=1}^k tr(W_j \rho_W  W_j^*) \log tr(H_j\rho_{\beta}  H_j^*)$$
\begin{equation}\label{p_var_denovo}
+\sum_{j=1}^k tr(W_j \rho_W  W_j^*) \log{\Big(\frac{(V_j\rho_{\beta} V_j^*)_{lm}}{(\rho_\beta)_{lm}}\Big)}
 \leq \log\beta,
\end{equation}
the {\bf basic inequality} associated to the potential $H\rho=\sum_i H_i\rho H_i^*$ and to the QIFS determined by $V_i$, $W_i$, $i=1,\dots, k$. Equality holds if for all $i,j,l,m$,
\begin{equation}
\frac{1}{\beta}tr(H_j\rho_{\beta}  H_j^*)\frac{(V_j\rho_{\beta} V_j^*)_{lm}}{(\rho_\beta)_{lm}}=\frac{tr(W_j V_i\rho_W V_i^* W_j^{*})}{tr(V_i\rho_W V_i^*)}
\end{equation}
\end{defi}

\qee

As before, $\rho_{\beta}$ is an eigenstate of $\mathcal{L}_H(\rho)$ and $\rho_W$ is the barycenter of the unique attractive, invariant measure for the Markov operator $\mathcal{V}$ associated to the QIFS $\mathcal{F}_W$. Given the classic inequality (\ref{desig_classica}) we want to compare it to the basic inequality (\ref{p_var_denovo}). More precisely, we would like to obtain operators $V_i$ that satisfy the following: given a matrix $A$ with positive entries and a stochastic matrix $Q$, there are  $H_i$ and $W_i$ such that inequality (\ref{p_var_denovo}) becomes inequality (\ref{desig_classica}). We have the following proposition.

\begin{pro}\label{prop_ba_cl}
Define
\begin{equation}
V_1=\left(
\begin{array}{cc}
1 & 0 \\
0 & 0
\end{array}
\right),\sp V_2=\left(
\begin{array}{cc}
0 & 1 \\
0 & 0
\end{array}
\right),\sp V_3=\left(
\begin{array}{cc}
0 & 0 \\
1 & 0
\end{array}
\right),\sp V_4=\left(
\begin{array}{cc}
0 & 0 \\
0 & 1
\end{array}
\right)
\end{equation}
Let $A=(a_{ij})$ be a matrix with positive entries and $Q=(q_{ij})$ a two-dimensional column-stochastic matrix. Define
\begin{equation}
H_{11}=\left(
\begin{array}{cc}
\sqrt{e^{a_{11}}} & \sqrt{e^{a_{11}}} \\
0 & 0
\end{array}
\right),\sp H_{12}=\left(
\begin{array}{cc}
\sqrt{e^{a_{12}}} & \sqrt{e^{a_{12}}} \\
0 & 0
\end{array}
\right)
\end{equation}
\begin{equation}
H_{21}=\left(
\begin{array}{cc}
0 & 0 \\
\sqrt{e^{a_{21}}} & \sqrt{e^{a_{21}}}
\end{array}
\right),\sp H_{22}=\left(
\begin{array}{cc}
0 & 0 \\
\sqrt{e^{a_{22}}} & \sqrt{e^{a_{22}}}
\end{array}
\right)
\end{equation}
and also
\begin{equation}
W_1=\left(
\begin{array}{cc}
\sqrt{q_{11}} & 0 \\
0 & 0
\end{array}
\right),\sp W_2=\left(
\begin{array}{cc}
0 & \sqrt{q_{12}} \\
0 & 0
\end{array}
\right)
\end{equation}
\begin{equation}
W_3=\left(
\begin{array}{cc}
0 & 0 \\
\sqrt{q_{21}} & 0
\end{array}
\right),\sp W_4=\left(
\begin{array}{cc}
0 & 0 \\
0 & \sqrt{q_{22}}
\end{array}
\right)
\end{equation}
Then the basic inequality associated to $W_i, V_i, H_i$, $i=1,\dots,4$, $l=m=1$ or $l=m=2$, is equivalent to the classic inequality associated to $A$ and $Q$.
\end{pro}

We use the following lemma.

\begin{lem}\label{lema_diag1}
For $V_i$ given by
\begin{equation}
V_1=\left(
\begin{array}{cc}
\sqrt{v_{11}} & 0 \\
0 & 0
\end{array}
\right),\sp V_2=\left(
\begin{array}{cc}
0 & \sqrt{v_{12}} \\
0 & 0
\end{array}
\right)
\end{equation}
\begin{equation}
V_3=\left(
\begin{array}{cc}
0 & 0 \\
\sqrt{v_{21}} & 0
\end{array}
\right),\sp V_4=\left(
\begin{array}{cc}
0 & 0 \\
0 & \sqrt{v_{22}}
\end{array}
\right)
\end{equation}
 where  $v_{ij}>0$, we have that the associated QIFS is such that $\rho_W$ and $\rho_\beta$ are diagonal density operators, for any choice of $W_i$ and $H_i$, $i=1,\dots, 4$.
\end{lem}
{\bf Proof of Lemma \ref{lema_diag1}} We have that $\rho_W$ is a fixed point of
$$\Lambda(\rho)=\sum_i tr(W_i\rho W_i^*)\frac{V_i\rho V_i^*}{tr(V_i\rho V_i^*)}$$
Writing
$$\rho=\left(
\begin{array}{cc}
\rho_{11} & \rho_{12} \\
\rho_{12} & \rho_{22}
\end{array}
\right),$$
we have that $\Lambda(\rho)=\rho$ leads us to
$$\frac{tr(W_1\rho W_1^*)}{tr(V_i\rho V_i^*)}\left(
\begin{array}{cc}
v_{11}\rho_{11} & 0 \\
0 & 0
\end{array}
\right)+\frac{tr(W_2\rho W_2^*)}{tr(V_2\rho V_2^*)}\left(
\begin{array}{cc}
v_{12}\rho_{22} & 0 \\
0 & 0
\end{array}
\right)$$
$$+\frac{tr(W_3\rho W_3^*)}{tr(V_3\rho V_3^*)}\left(
\begin{array}{cc}
0 & 0 \\
0 & v_{21}\rho_{11}
\end{array}
\right)+\frac{tr(W_4\rho W_4^*)}{tr(V_4\rho V_4^*)}\left(
\begin{array}{cc}
0 & 0 \\
0 & v_{22}\rho_{22}
\end{array}
\right)=\left(
\begin{array}{cc}
\rho_{11} & \rho_{12} \\
\rho_{12} & \rho_{22}
\end{array}
\right)$$
Then $\rho_{12}=0$ and so $\rho_W$ is diagonal. In a similar way we prove $\rho_\beta$ is diagonal.
\qed

{\bf Proof of Proposition \ref{prop_ba_cl}} Let $V_i$, $W_i$, $i=1,\dots, 4$ and $H_{ij}$, $i,j=1,2$ as in the statement of the proposition. A simple calculation shows that
\begin{equation}
tr(H_{ij}\rho_\beta H_{ij}^*)=e^{a_{ij}}
\end{equation}
(since $\rho_\beta$ is diagonal, by lemma \ref{lema_diag1}). By example \ref{caso_util_a1}, the choice of $V_i$ and $W_i$ we made is such that the entropy $h_V(W)$ reduces to the Markov chain entropy. Then a calculation yields
\begin{equation}
\frac{(V_i\rho_\beta V_i^*)_{11}}{(\rho_\beta)_{11}}=\frac{(\rho_\beta)_{11}}{(\rho_\beta)_{11}}=1
\end{equation}
In a similar way,
\begin{equation}
\frac{(V_i\rho_\beta V_i^*)_{22}}{(\rho_\beta)_{22}}=\frac{(\rho_\beta)_{22}}{(\rho_\beta)_{22}}=1
\end{equation}
Then from the basic inequality with  $l=m=1$ or $l=m=2$ we get
\begin{equation}\label{pv_pr01bb}
h_V(W)+\sum_j tr(W_j\rho_W W_j^*)\sum_i \frac{tr(W_i V_j\rho_W V_j^* W_i^{*})}{tr(V_j\rho_W V_j^*)}\log{tr(H_i\rho_{\beta} H_i^*)}\leq \log\beta
\end{equation}
Finally, since $tr(H_{ij}\rho_{\beta} H_{ij}^*)=e^{a_{ij}}$ and $Q\pi=\pi$, we conclude that (\ref{pv_pr01bb}) becomes (\ref{desig_classica}).

\qed

\begin{examp}
Let
$$H_1=\left(
\begin{array}{cc}
2i & 2i \\
0 & 0
\end{array}
\right), \sp H_2=I,\sp H_3=\left(
\begin{array}{cc}
i\sqrt{2} & i\sqrt{2} \\
0 & 0
\end{array}
\right),\sp H_4=I$$
Then
$$H_1^*=\left(
\begin{array}{cc}
-2i & 0 \\
-2i & 0
\end{array}
\right),\sp H_2^*=I,\sp H_3^*=\left(
\begin{array}{cc}
-i\sqrt{2} & 0 \\
-i\sqrt{2} & 0
\end{array}
\right),\sp H_4^*=I$$
If we suppose the $V_i$ are the same as from proposition \ref{prop_ba_cl}, we have that $\rho_\beta$ is diagonal, so
$$tr(H_1\rho_\beta H_1^*)=4 ,\sp tr(H_2\rho_\beta H_2^*)=1,\sp tr(H_3\rho_\beta H_3^*)=2 ,\sp tr(H_4\rho_\beta H_4^*)=1$$
Then $\mathcal{L}_H(\rho)=\beta\rho$ leads us to
$$4\rho_{11}+\rho_{22}=\beta\rho_{11}$$
$$2\rho_{11}+\rho_{22}=\beta\rho_{22}$$
A simple calculation gives
$$\beta=\frac{5+\sqrt{17}}{2}$$
with eigenstate
$$\rho_\beta=\frac{4}{7+\sqrt{17}}\left(
\begin{array}{cc}
\frac{3+\sqrt{17}}{4} & 0 \\
0 & 1
\end{array}
\right)$$
\end{examp}

\qee

We want to calculate the $W_i$ which maximize the basic inequality (\ref{p_var_denovo}). Recall that from proposition \ref{prop_ba_cl}, the choice of $V_i$ we made is such that
$$\frac{(V_j\rho_{\beta} V_j^*)_{lm}}{(\rho_\beta)_{lm}}=1,$$
So
\begin{equation}
h_V(W)+ \sum_{j=1}^k tr(W_j \rho_W  W_j^*) \log tr(H_j\rho_{\beta}  H_j^*)\leq \log\beta
\end{equation}
and equality holds if and only if, for all $i,j,l,m$,
\begin{equation}\label{a_cond_pw}
\frac{1}{\beta}tr(H_j\rho_{\beta}  H_j^*)\frac{(V_j\rho_{\beta} V_j^*)_{lm}}{(\rho_\beta)_{lm}}=\frac{tr(W_j V_i\rho_W V_i^* W_j^{*})}{tr(V_i\rho_W V_i^*)}
\end{equation}
Choose, for instance, $l=m=1$. Then condition (\ref{a_cond_pw}) becomes
\begin{equation}
\frac{1}{\beta}tr(H_j\rho_{\beta}  H_j^*)=\frac{tr(W_j V_i\rho_W V_i^* W_j^{*})}{tr(V_i\rho_W V_i^*)}
\end{equation}
To simplify calculations, write $\widehat{W}_i=W_i^*W_i$ and $\widehat{W}_i=(w_{ij}^i)$. Then we get
\begin{equation}
\frac{tr(H_i\rho_\beta H_i^*)}{\beta}=w_{11}^i=w_{22}^i,\sp i=1,\dots, 4
\end{equation}
So we conclude
\begin{equation}
W_i=\frac{1}{\sqrt{\beta}}\left(
\begin{array}{cc}
\sqrt{tr(H_i\rho_\beta H_i^*)} & 0 \\
0 & \sqrt{tr(H_i\rho_\beta H_i^*)}
\end{array}
\right),\sp i=1,\dots, 4
\end{equation}
That is,
\begin{equation}
W_1=\frac{2}{\sqrt{\beta}}I,\sp W_2=\frac{1}{\sqrt{\beta}}I,\sp
W_3=\frac{\sqrt{2}}{\sqrt{\beta}}I,\sp
W_4=\frac{1}{\sqrt{\beta}}I
\end{equation}
Note that
$$\sum_i W_i^*W_i=\frac{4+\sqrt{2}}{\sqrt{\beta}}I\neq I$$
To solve that, we renormalize the potential. Define
\begin{equation}
\tilde{H}_i:=\sqrt{\alpha} H_i,\sp \alpha:=\frac{\sqrt{\beta}}{4+\sqrt{2}}
\end{equation}
Then a calculation shows that $\mathcal{L}_{\tilde{H}}(\rho)=\tilde{\beta}\rho$ gives us the same eigenstate as before, that is $\rho_{\tilde{\beta}}=\rho_\beta$. But note that the associated eigenvalue becomes $\tilde{\beta}=\alpha\beta$. Now, note that it is possible to renormalize the $W_i$ in such a way that we obtain $\tilde{W}_i$ with $\sum_i \tilde{W}_i^*\tilde{W}_i=I$, and that these maximize the basic inequality for the $H_i$ initially fixed. In fact, given the renormalized $\tilde{H}_i$, define
\begin{equation}
\tilde{W}_i=\sqrt{\alpha} W_i,\sp i=1,\dots, 4
\end{equation}
Note that $\sum_i \tilde{W}_i^*\tilde{W}_i=I$. Also we obtain
\begin{equation}
h_V(\tilde{W})+ \sum_{j=1}^k tr(\tilde{W}_j \rho_{\tilde{W}}  \tilde{W}_j^*) \log tr(\sqrt{\alpha}H_j\rho_{\beta} \sqrt{\alpha}H_j^*)\leq \log\alpha\beta
\end{equation}
which is equivalent to
\begin{equation}
h_V(\tilde{W})+ \sum_{j=1}^k tr(\tilde{W}_j \rho_{\tilde{W}}  \tilde{W}_j^*) \log (\alpha tr(H_j\rho_{\beta} H_j^*))\leq \log\alpha+\log\beta
\end{equation}
That is
$$
h_V(\tilde{W})+ \sum_{j=1}^k tr(\tilde{W}_j \rho_{\tilde{W}}  \tilde{W}_j^*) \log \alpha
$$
\begin{equation}
+\sum_{j=1}^k tr(\tilde{W}_j \rho_{\tilde{W}}  \tilde{W}_j^*) \log tr(H_j\rho_{\beta} H_j^*)\leq \log\alpha+\log\beta,
\end{equation}
and if we cancel $\log\alpha$ on both sides, we get the same inequality as for the nonrenormalized $H_i$. As we have seen before, such $\tilde{W}_i$ gives us equality. Hence
\begin{equation}
h_V(\tilde{W})+ \sum_{j=1}^k tr(\tilde{W}_j \rho_{\tilde{W}}  \tilde{W}_j^*) \log tr(H_j\rho_{\beta} H_j^*)=\log\beta
\end{equation}
\qee

\section{Remarks on the problem of pressure and quantum mechanics}\label{sec_analysis4}

One of the questions we are interested in is to understand how to formulate a variational principle for pressure in the context of quantum information theory. An appropriate combination of such theories could have as a starting point a relation between the inequality for positive numbers
$$-\sum_i q_i\log q_i+\sum_i q_i\log p_i\leq 0,$$
(lemma  \ref{lemalog2}, seen in certain proofs of the variational principle of pressure), and the QIFS entropy. We have carried out such a plan and then we have obtained the basic inequality, which can be written as
\begin{equation}\label{pvar_ver1yy}
h_V(W)+ \sum_{j=1}^k \log \Big( tr(H_j\rho_{\beta}  H_j^*)tr(V_j\rho_{\beta} V_j^*)\Big) tr(W_j \rho_W  W_j^*)
 \leq \log\beta
\end{equation}
where equality holds if and only if for all $i,j$,
\begin{equation}\label{cond_igualdadeyy}
\frac{1}{\beta}tr(H_j\rho_{\beta}  H_j^*)tr(V_j\rho_{\beta}  V_j^*)=\frac{tr(W_j V_i\rho_W V_i^* W_j^{*})}{tr(V_i\rho_W V_i^*)}
\end{equation}
As we have discussed before, it is not clear that given any
dynamics, we can obtain a measure such that we can reach the
maximum value $\log\beta$. Considering  particular cases we can
suppose, for instance, that the $V_i$ are unitary. In this way we
combine in a natural way a problem of classic thermodynamics,
with an evolution which has a quantum character. In this
particular setting, we have for each $i$ that
$V_iV_i^*=V_i^*V_i=I$ and then the basic inequality becomes
\begin{equation}\label{pvar_ver1ii}
h_V(W)+ \sum_{j=1}^k tr(W_j \rho_W  W_j^*)\log tr(H_j\rho_{\beta}  H_j^*)
 \leq \log\beta
\end{equation}
and equality holds if and only if for all $i,j$,
\begin{equation}\label{cond_igualdadeii}
\frac{1}{\beta}tr(H_j\rho_{\beta}  H_j^*)=tr(W_j V_i\rho_W V_i^* W_j^{*})
\end{equation}
We have the following:
\begin{lem}
Given a QIFS with a unitary dynamics (i.e., $V_i$ is unitary for each $i$),
there are $\hat{W}_i$ which maximize (\ref{pvar_ver1yy}), i.e., such that
\begin{equation}\label{pvar_ver1iij}
h_V(\hat{W})+ \sum_{j=1}^k tr(\hat{W}_j \rho_{\hat{W}}  \hat{W}_j^*)\log tr(H_j\rho_{\beta}  H_j^*)
= \log\beta
\end{equation}
\end{lem}
{\bf Proof} Define, for each $j$,
\begin{equation}
\hat{W}_j:=\sqrt{\frac{1}{\beta}tr(H_j\rho_{\beta}H_j^*)}I
\end{equation}
where $I$ is the identity. The equality condition (\ref{cond_igualdadeii}) is satisfied by such $\hat{W}_j$, so the lemma follows.
\qed

{\bf Remark} The above lemma also holds for the basic inequality in coordinates, given by (\ref{p_var_denovo}). Also, it is immediate to obtain a similar version of the above lemma for any QIFS such that the $V_i$ are multiples of the identity, and also for QIFS such that $\rho_W$ fixes each branch of the QIFS, that is, satisfying, for each $i$,
$$\frac{V_i\rho_W V_i^*}{tr(V_i\rho_W V_i^*)}=\rho_W$$

\qee

\section{Concluding remarks}

Considering the QIFS setting, we defined a concept of entropy and a Ruelle operator in such a way that we are able to get some analogous results to the classical Thermodynamic Formalism. Such Ruelle operator admits a positive eigenvalue, which gives us an upper bound for the pressure (entropy plus a potential) associated to the QIFS. Our configuration space is the set of density matrices. We did not consider the usual space of symbols or a shift operator, as it is assumed in the Ruelle-Perron-Frobenius theory. We have replaced the dynamics given by the shift with the one given by the inverse branches of the iterated functions (which are defined by a set of operators).

The references \cite{lozinski} and \cite{wsbook} are of fundamental importance in
our investigation.

A starting point for further investigation could be to study more properties of the QIFS entropy, such as convexity and subadditivity. Also, a natural question is to ask whether it is possible to consider a QIFS  acting in an infinite tensor product of finite Hilbert spaces  which would be the analogous of considering the full Bernoulli space.

In a forthcoming paper we are going to consider relative entropies and quantum conditional expectations.

\end{document}